\documentclass[12pt]{article}
\usepackage{amsmath,amssymb,latexsym}
\usepackage{graphicx}
\usepackage{epstopdf}

\def\qed{$\Box$}

\def\pf{{\bf Proof }}

\begin{document}
\title{Parametric argument principle and its applications to CR functions and manifolds}
\author{Mark L.~Agranovsky}
\maketitle
\begin{center}
Bar-Ilan University
\end{center}

\newcommand{\bt}{\begin{Theorem}}
\newcommand{\et}{\end{Theorem}}
\newcommand{\bi}{\begin{itemize}}
\newcommand{\ei}{\end{itemize}}
\newcommand{\bea}{\begin{eqnarray}}
\newcommand{\eea}{\end{eqnarray}}
\newtheorem{Definition}{Definition}[section]
\newtheorem{Theorem}[Definition]{Theorem}
\newtheorem{Lemma}[Definition]{Lemma}
\newtheorem{Exercise}{\sc Exercise}[section]
\newtheorem{Proposition}[Definition]{Proposition}
\newtheorem{Corollary}[Definition]{Corollary}
\newtheorem{Problem}[Definition]{Problem}
\newtheorem{Example}[Definition]{Example}
\newtheorem{Remark}[Definition]{Remark}
\newtheorem{Remarks}[Definition]{Remarks}
\newtheorem{Question}[Definition]{QuesMotition}
\newtheorem{Statement}[Definition]{Statement}
\newcommand{\be}{\begin{equation}}
\newcommand{\ee}{\end{equation}}
\def\pf{{\bf Proof }}
\vskip.25in

\setcounter{equation}{0}

\begin{abstract}
It follows from the classical argument principle that if a holomorphic mapping $f$ from the unit disc $\Delta$ to $\mathbb C$ or, more generally, to $\mathbb C^n,$ smooth in the closed disc,
is homologicaly trivial on the unit circle  (i.e.,  $H_1(\gamma)=0, \gamma=f(S^1),$ which is equivalent to either $\gamma$ being a point or $\gamma$ having endpoints, i.e., $\partial \gamma \neq \emptyset$), then $f=const$, in other words,  the image of the unit disc degenerates to a point.
We establish a parametric version of this fact, for a variety of holomorphic mappings from $\Delta$ to $\mathbb C^n$ in place of  a single mapping. We find conditions for
a holomorphic mapping of the unit disc, depending on additional real parameters, under which
homological triviality of the boundary image  implies collapse of the dimension of the image of the interior. As an application, we obtain
estimates of dimensions of complex tangent bundles of real submanifolds in $\mathbb C^n$, in terms of zero moment conditions on families of closed curves covering the manifold.
Applying this result to the graphs of functions, we obtain  solution of several known problems about characterization of holomorphic CR functions in terms of moment conditions on families of curves.
\end{abstract}



\section {Content}

\begin{enumerate}
\item Introduction
\item Content of the article
\item
Basic definitions and notations
\item
Formulation of Main Theorem
\item
Proof of Theorem \ref{T:Main}
\item Theorem \ref{T:Main} in terms of attached analytic discs: detecting positive CR dimensions
\item Lower bounds for higher CR dimensions
\item Characterization of holomorphic manifolds and their boundaries
\item Morera type theorems for CR functions
\item Concluding remarks
\end{enumerate}

\section{Introduction} \label{S:intro}
This article is devoted to study holomorphy by means of integral-geometric tools. Briefly, the problem which is studied here can be formulated as follows:
given a real submanifold of a complex manifold, is it possible to detect a complex structure on the manifold in terms of zero moment conditions on a family of closed curves covering the manifold?

The starting point for this article were the works of the author  \cite{A1, A2}. There,
two open problems, mentioned in the abstract, were solved in real-analytic case. These problems concern characterization of holomorphic functions in a planar domain and characterization of boundary values of holomorphic functions in a domain in $\mathbb C^n,$  in terms of complex moment conditions on closed curves.

In this article we show that both problems, after being reformulated in terms of the graphs of the functions, become just two special cases of a more general problem of deriving lower bounds for CR dimensions of real manifolds in $\mathbb C^n$ from the moment conditions on families of closed curves.

Solution of the latter problem appeared to be intimately related to a generalization of the classical argument principle, from a single holomorphic mappings to  varieties of those.
We prove this generalization, which we  call  parametric argument principle, and then show how it leads to solutions of the above mentioned problems about CR functions and manifolds, and of more general results.

\subsection{Motivation and description of Main Theorem}\label{S:Motivation}

We will be denoting $\Delta:=\{z \in \mathbb C:|z|<1\}$-the unit complex disc and $\partial \Delta=S^1=\{z \in \mathbb C:|z|=1\}$ - the unit circle.

Let $\Phi$ be an analytic function in the open unit disc $\Delta \subset \mathbb C,$
of the class $C^1$ in the closed disc $\overline \Delta.$ The classical argument principle states that for any complex number
$b \notin \Phi (\partial \Delta)$ the following identity  holds
$$
N_{\Phi}(b)=W(\Phi(\partial \Delta),b),
$$
where $N_{\Phi}(b)$ denotes the number of zeros, counting multiplicities,  of the function $\Phi(z)-b$ in $\Delta$
and $W(\gamma,b)$ stands for the winding (rotation) number of a curve $\gamma$ with respect to the point $b  \notin \gamma.$
The winding number $W(\Phi(\partial \Delta), b)$ is defined as the change of argument of $\Phi(z)$ along the unit circle $\partial \Delta,$ divided by $2\pi$,  and can  be evaluated via logarithmic residue:
$$
W(\gamma,b)=\frac{1}{2\pi i}\int\limits_{\partial \Delta}\frac{d\Phi(z)}{\Phi(z)-b}.
$$

Notice that if $\Phi$ is holomorphic in the closed unit disc then $\Phi-b$ is allowed to have (isolated) zeros on the unit circle $\partial\Delta$. The multiplicities of the boundary zeros contribute to $N_{\Phi}(b)$
with the factor $\frac{1}{2}.$

A generalization of the argument principle to mappings to $\mathbb C^n$ addresses to linking numbers, or linking coefficients (cf. \cite{AW}), which generalize the notion of winding number. By definition, the linking number of a closed curve $\gamma$  and a $2n-2$-manifold $V,$ disjoint from $\gamma$, equals the sum of intersection indices of $V$ and $A,$ and equals zero if $V \cap A=
\emptyset,$ where $A$ is a 2-chain bounded by $\gamma.$ The fundamental fact is that if $A$ and $V$ are holomorphic, then the intersection indices are positive (cf.\cite{GH}, p.62) and
therefore the linking number equals to the number of intersection points, counting multiplicities:
$$link(\gamma,V)=\# (A \cap V)$$
Also,  if $V$ is given as the null variety of a holomorphic function, $V=P^{-1}(0),$ then the linking number can be computed via the logarithmic residue formula: $$link(\gamma,V)=\frac{1}{2\pi}\int\limits_{\gamma}\frac{dP}{P}.$$

The facts that the intersection indices are positive and, correspondingly, the linking numbers coincide with the number of intersection points,
can be regarded as a general form of the argument principle.
We are interested in the following corollary, which can be viewed as a weak form of the argument principle.

\begin{Proposition}(Corollary of Argument Principle)\label{P:argument_principle}
Let $\Phi:\overline \Delta \to \mathbb C^n$ be $C^1$-mapping, holomorphic in $\Delta.$
If  $H_1(\Phi(\partial \Delta))=0$, then $\Phi^{\prime}(z)=0$, i.e. $\Phi=const.$
\end{Proposition}

It is instructive to see  the proof of that fact.

\pf
First, we prove that the images of the unit circle and of the unit disc coincide: $\Phi(\overline \Delta)=\Phi(\partial\Delta).$ This will imply that $\Phi=const.$ Indeed, otherwise $\Phi^{\prime}(\zeta) \neq (0,\cdots,0)$
for some $\zeta \in \Delta$ and then the image of $\Phi,$ near $\Phi(\zeta),$ is a two-dimensional manifold, while $\Phi(\partial\Delta)$ is at most one-dimensional.

If, on the contrary,  $\Phi(\overline \Delta) \neq \Phi(\partial\Delta),$
then there exists a point $b \in \Phi(\overline\Delta)\setminus \Phi(\partial\Delta).$  Since smooth curves in $\mathbb C^n$ are rationally convex (e.g., \cite{Gam}, Theorem 2.7),
there exists a polynomial $P$ in $\mathbb C^n$ such that $P(b)=0$ but $P(z) \neq 0$  for $z \in \Phi(\partial \Delta).$

Since the image $\Phi(\partial \Delta)$  is homologically trivial, then  $link(\Phi(\partial\Delta),V)=0,$ where $ \ V=P^{-1}(0).$ On the other hand, the intersection
$\Phi(\Delta) \cap V$ is nonempty, because it contains the point $b$, and therefore $link(\Phi(\partial\Delta),V) >0,$ due to the argument principle. This contradiction completes the proof.
\qed

In fact, the crucial condition for the mapping $\Phi$ is $\deg \ \Phi\vert_{S^1}=0.$ Here $\deg$ means the topological degree of the mapping $\Phi:S^1 \to \Phi(S^1).$
Hence the above statement sounds as follows:
{\it if $\Phi \in C^1(\overline \Delta)$ is holomorphic in $\Delta$ and $\deg \Phi\vert_{S^1}=0,$ then $\Phi= const.$}

Thus,  the classical Argument Principle implies that if a holomorphic mapping of the unit disc, smooth in its closure, homologically degenerates on the boundary (the topological degree is zero), then it dimensionally degenerates on the interior, meaning that
$rank \ d\Phi(\zeta)=0 < dim \Delta=2, \ \zeta \in \Delta,$ and, correspondingly, the image collapses.

The main goal of this article is  to generalize this phenomenon from a single mapping to families of holomorphic mappings. Namely, we consider a variety
$\overline\Delta$ to $\mathbb C^n,$ real-analytically depending on the parameter $t.$ This parameter runs over a compact real-analytic connected oriented $k-$dimensional  manifold $M.$

For the sake of simplicity of the exposition, we resctrict ourselves here by the case when the parametrizing manifold $M$ has no boundary (is closed) . 

The family $\Phi_t$ defines a smooth mapping $\Phi(\zeta,t)=\Phi_t(\zeta)$ from $\overline \Delta \times M$ to $\mathbb C^n$, holomorphic in the variable $\zeta \in \Delta.$ We will be assuming that the target space has enough room, namely, that   $2n \geq \dim (\Delta \times M)=dim M+2.$

We want to understand when the homological degeneracy
of the mapping $\Phi$ on the boundary $\partial \Delta \times M$ implies dimensional degeneracy of $\Phi$ in the interior $\Delta \times M.$
In general terms, the main theorem (Theorem \ref{T:Main}) of this article can be formulated  as the following

\noindent
{\bf Parametric Argument Principle:} {\it If the mapping $\Phi(\zeta,t),$ holomorphic in $\zeta,$ homologically
degenerates on the boundary $\partial \Delta \times M$ but does not homologically degenerate in the interior $\Delta \times M$, then $\Phi$ dimensionally degenerates on the whole manifold $\overline \Delta \times M.$}

Let us explain the terminology and comment on the statement.
By homological degeneracy on the boundary is meant the following :

\noindent
{\bf Condition 1.} The induced homomorphism $(\Phi\vert_{\partial \Delta \times M})_* : H_{k+1}(\partial\Delta \times M) \to H_{k+1}(\Phi(\partial \Delta \times M)), \ k=dim M,$
is trivial: $(\Phi\vert_{\partial \Delta \times M})_*=0.$

By the homological non-degeneracy in the interior is meant the following:

\noindent
{\bf Condition 2.} The induced homomorphism $\Phi_*:H_{k}(\Delta \times M) \to H_k(\Phi(\Delta \times M))$ is nontrivial: $\Phi_* \neq 0.$ Equivalently, this means that the orbit
$\Phi (\{\zeta_0\} \times M=\{\Phi(\zeta_0,t): \ t \in M \},$ of a point $\zeta_0 \in \overline \Delta$, when the parameter $t$ runs the $k-$ dimensional manifold $M,$ represents a nonzero class in
the $k-$ dimensional ($k=dim M$) homology group $H_k (\Phi(\Delta \times M))$ of the image.  Clearly, if the above condition holds for some $\zeta_0,$ then it holds for all $\zeta_0 \in \overline \Delta.$

\begin{Remark} In order not to make the presentation more complicated, we consider in this article only the case of closed parametrizing manifolds $M$. The case of manifolds $M$ with nonempty boundaries requires
considering relative homology groups.
\end{Remark}

Under Conditions 1 and 2, the conclusion is that the rank of the mapping $\Phi$ degenerates, i.e., each point in $\overline\Delta \times M$ is critical:
$rank \ d\Phi(\zeta,t) <  k+2=\dim \Delta \times M, \ (\zeta,t) \in \overline \Delta \times M.$

\noindent
{\bf Heuristics.} The phenomenon under discussion might be explained as follows.

Condition 1 is the same as in the classical parameter-free version, namely, it requires that the holomorphic mapping $\Phi$ is homologically trivial on the boundary of the unit disc $\Delta.$

Condition 2 rules out the parameter factor. Namely,  Condition 2 of the homological non-degeneracy in the interior excludes the possibility that the boundary homological degeneracy  occurs due to the parameter $t$, and results in becoming the holomorphic variable $\zeta$ the main player . This triggers the  argument principle for holomorphic mappings and consequently leads to falling down the rank of the mapping and, correspondingly, to the collapse of the image, as in the parameter-free case.

Let us illustrate this idea one the following model situation.

{\noindent}
{\bf A model example.}

Let $\Phi(\zeta,t)=(\Phi_1(\zeta,t),\Phi_2(\zeta,t))$ be a mapping of the solid torus $\overline \Delta \times S^1, \ S^1=\{t \in \mathbb C:|t|=1,\}$ into itself such that
$\Phi(\zeta,t)$ is holomorphic with respect to $\zeta \in \overline \Delta.$

Condition 2 says that for $\zeta_0 \in \partial \Delta=S^1,$ the orbit $\{\Phi(\zeta_0, t)=(\Phi_1(\zeta_0,t),\Phi_2(t)): t \in S^1\}$  is a homologically nontrivial
loop in  $\Phi(\overline D \times S^1).$
On the other hand, Condition 1 requires that the image $\Phi(T^2)$ of the boundary torus represents the zero homology class in $H_2(\Phi(T^2).$  Since the curves $C_{\zeta_0}=\Phi(\{\zeta_0\}\times S^1)$  are  homologically nontrivial (Condition 2),
the second generator, the curves
$C^t=\Phi(S^1 \times \{t\})$ must be homologically trivial loops in $\mathbb C.$ But then, by the Argument Principle, the holomorphic mappings $\Phi_t(\zeta)=\Phi(\zeta,t)$ are constant with respect to  $\zeta.$ Therefore, $\Phi(\zeta,t)=(\Phi_1(t),\Phi_2(t))$
and hence the image collapses, as it is at most one-dimensional.

Simple examples show that Condition 2 is crucial and can not be dropped. For instance, take  $M=S^1$ realized as the unit circle $t_1^2+t_2^2=1$ in $\mathbb R^2.$  Define $\Phi(\zeta,t)=(\zeta, t_1).$ Then $k=1,$ the boundary image is $\Lambda=\Phi(\partial \Delta \times S^1)=S^1 \times [-1,1]$ and the induced homomorphism $\Phi_*:H_2(\partial \Delta \times M) \to H_2(\Lambda)=0$ is trivial, so Condition 1 holds. On the other hand,  Condition 2 is not fullfiled,
because $H_1 (\Phi(\Delta \times M))=H_1(\Delta \times [-1,1])=0$ and hence the induced homomorphism of the first homology groups is $\Phi_*=0$. The conclusion also fails since the image $\Phi(\Delta \times M)=\Delta \times [-1,1]$ is three dimensional, $rank \Phi=3,$  so that  $\Phi$ does not degenerate in the sense of the rank.

\subsection{Applications}\label{S:Intro-appl}
The geometric applications of Theorem \ref{T:Main} are concerned with characterization of the boundary image $\Lambda:=\Phi(\partial \Delta \times M)$ which is assumed a real manifold in $\mathbb C^n.$

First of all, for each fixed $t$ the image of the unit disc $\Delta$ is an analytic disc $D_t:=\Phi(\Delta \times \{t\})$ in $\mathbb C^n.$ This disc is {\it attached} to $\Lambda, $ meaning that
$\gamma_t:=\partial D_t \subset \Lambda.$ Moreover, the manifold $\Lambda$ is fully covered by the closed curves $\gamma_t.$

The key point for applications is  that the conclusion of Theorem \ref{T:Main} about the degeneracy of the differential $d\Phi$  means
 $\mathbb C-$ linear dependence in the tangent spaces to $\Lambda.$ In turn, this implies that the  tangent spaces of $\Lambda$ contain nontrivial complex subspaces and hence $\Lambda$ possesses a partial holomorphic (CR) structure.

Thus,  speaking in geometric terms, Theorem \ref{T:Main}  characterizes  such varieties of analytic discs which induce CR, i.e., partially complex,  structure  
on a real manifold in $\mathbb C^n,$ swept by the boundaries of the  discs.
Moreover, it appears possible to estimate the dimension of the above CR structure (CR dimension) in terms of the properties of the varieties of attached analytic discs.

An alternative, analytic,  interpretation of the result is as follows. Since $\int_{\gamma}\omega=0$ for any  closed curve, that bounds an analytic disc, and for each holomorphic 1-form $\omega,$
the above application can be regarded as a {\it version of Morera theorem for manifolds}. Namely,  it detects holomorphic structure on $\Lambda$ from complex moment conditions $\int_{\gamma}\omega=0$ on a family of
closed curves $\gamma$ covering $\Lambda.$

The discussed result about estimating CR dimension leads to solutions, in real-analytic category, of number or problems which turned out special cases of the above mentioned Morera-type theorem for manifolds. In fact, just attempts
to solve those problem had led the author to Parametric Argument Principle which appears to be in a core of the problems.
Let us describe those applications.

\bigskip

\noindent
{\bf Strip-problem.}

The first problem, called strip-problem,  concerns Morera-like characterization of analytic functions of one complex variables.
In the general setting it sounds as follows.

Let $\gamma_t$ be a 1-parameter family of Jordan curves in the complex plane. Let
$f$ be a continuous function on $\Omega=\cup_t \gamma_t$. Suppose that all complex  moments on the curves $\gamma_t$ vanish:
$$\int\limits_{\gamma_t}f(z)z^kdz=0, \ k=0,1,\cdots$$
for all values of the parameter $t.$  Must $f$ be holomorphic in $\Omega?$

The moment condition is known to be equivalent to existence of holomorphic extension of the function $f$ from any curve $\gamma_t$ to the Jordan domain $\Omega_t$ bounded by $\gamma_t.$
The  main question is when these extensions agree on intersections and thus  $f$ is holomorphic on $\Omega?$
The name "strip-problem" comes from the original formulation of the problem, where $\gamma_t$ were the circles inscribed in a given strip. Clearly, not any family of Jordan curves test holomorphicity in the
above sense. For example, the family of circles $\{|z|=t\}$  does not, since the function $f(z)=|z|$ extends holomorphically from each such circle (as a constant) but is not holomorphic. So, the problem is to characterize those families which do test holomorphic functions.

The strip-problem has quite long history and bibliography (see the survey \cite{AZ}). In the articles  \cite{AV},
\cite{AG},\
\cite{E2},\cite{E1},
\cite{G1}, \cite{G3}, \cite{G4},\cite{G6},\cite{G7},
\cite{T1},\cite{T2}, \cite{Z1} essentially two situations have been considered and solved till recently: 1) the case of rotation invariant families of curves
(this assumption makes methods of Fourier series applicable), 2) the case of  circles.

However, even for group-invariant families of circles, in the case of non-compact groups (e.g.,translations along a line) the problem was open for a long time. Relatively recently, the problem has been pushed forward due to reformulating it as a problem in $\mathbb C^2$. This approach was
used in \cite{AG} to solve the problem for rational, with respect to the real variables, functions and for generic families of circles. Then bringing tools of CR theory led
A. Tumanov to  a complete solution for the case of families of circles  \cite{T1,T2}.  A crucial condition on a family of circles was that they surround no common point.

In full generality, for generic families of  Jordan curves, though under assumptions of real analyticity, the strip-problem was solved by the author in
\cite{A1,A2}. It was proved that for any regular real-analytic family of Jordan curves, surrounding {\it no common point}, the strip-problem has positive answer. The proof was based on topological ideas related to the argument principle, so that the concept of parametric argument principle was implicitly present in \cite{A1,A2}.
It should be mentioned that certain hints to topological nature of the problem appeared already in Tumanov's paper \cite{T2} where such tools as
indices of the curves were exploited. We present here the results of \cite{A1,A2}, as well as their  generalizations to $\mathbb C^n$ (Corollary  \ref{C:4.2}, Corollary \ref{C:n-dim-strip}), as a corollary of the parametric argument principle, proved in this article.

\bigskip

\noindent
{\bf Globevnik-Stout conjecture}

The second problem that served as a motivation for this article is Globevnik-Stout conjecture \cite{GS}.
It is similar to the strip-problem,  but is formulated in higher dimensions and
concerns testing boundary values of holomorphic functions, or more generally, CR functions.

The question goes back to the work of of Nagel and Rudin \cite{NR} where they proved that a continuous function $f$ on the unit sphere in $\mathbb C^n$ extends holomorphically in the unit ball $B^n$  if $f$ possesses holomorphic extension into each complex line tangent to a fixed ball contained in $B^n$. The conjecture in \cite{GS} was that the result generalizes for the case of two arbitrary  enclosed domains $D^{\prime} \subset D \subset \mathbb C^n$.
Partial results were obtained for open families of complex lines \cite{GS}, \cite{AS1}, or for complex geodesics in Kobayashi metric \cite{BTZ}.

In \cite{A1,A2} the Globevnik-Stout conjecture was proved in real-analytic category, even in a much broader formulation. Namely, the conjecture  was proved for general families of attached analytic discs, of which the family of complex lines tangent to a given surface is a particular case. In this article, we present a proof of even more general statement  (Theorem \ref{T:G_S}), as a corollary of the results we are going to describe next.

\bigskip

\noindent
{\bf One-dimensional Holomorphic Extension Problem}

Let us formulate a more general problem which contains
both the strip-problem and Globevnik-Stout conjecture as special cases.
\begin{Definition} Let $\Omega$ be a real smooth manifold in $\mathbb C^n.$ Let $\mathcal F=\{D\}$ be a family of analytic discs $D$ attached
to $\Omega$ so that $\Omega=\cup_{D \in \mathcal F}\partial D.$ We say that a function $f \in C(\Omega)$ satisfies the property of one-dimensional
holomorphic extension with respect to $\mathcal F$ if for any $D \in \mathcal F$ there exists a function $F_D \in C(\overline D)$, holomorphic in $D,$ such that $f(z)=F_D(z)$ for $z \in \partial D.$
\end{Definition}

\noindent
{\bf Problem.} {\it Let $\Omega \subset \mathbb C^n$ be a CR manifold. Characterize the families $\mathcal F$ of analytic discs attached to $\Omega$
such that if $f \in C^1(\Omega)$ satisfies the  property of one-dimensional holomorphic extension with respect to $\mathcal F,$ then $f$ is a  CR function on $\Omega.$ }

Remind that a real manifold $\Omega$ in $\mathbb C^n$ is CR if the complex tangent bundle has a constant dimension and a smooth function $f$ on $\Omega$ is CR if $f$ satisfies tangential
$\overline\partial-$ conditions in complex directions.

The strip-problem and Globevnik-Stout conjecture are just special cases of One-dimensional Holomorphic Extension Problem.
Indeed, the strip-problem is the particular case  when the dimension $n=1.$
In this situation $\Omega$ is the domain covered by the curves $\gamma_t$. The family $\mathcal F$ consists of
the Jordan domains $D_t \subset \mathbb C$ bounded by the curves $\gamma_t,$ i.e., in this case the analytic discs attached to $\Omega$  are simply contained in $\Omega.$
The property of one-dimensional holomorphic extension with respect to $\mathcal F$
is  exactly the property of analytic extendibility to the interiors of the curves $\gamma_t$ and  CR-functions in the question are just holomorphic functions in $\Omega.$

In Globevnik-Stout conjecture,  the manifold $\Omega=\partial D$ is the boundary of a domain $D$ in $\mathbb C^n.$ The family $\mathcal F$ of attached analytic discs consists of the intersections $D_L=L \cap D$ of the larger domain $D$ and the complex lines $L$ tangent to the smaller domain $D^{\prime}.$ The condition  is that $f$ satisfies the one-dimensional extension property with respect to $\mathcal F,$ and the conclusion that
$f$ is a CR function. This is equivalent, due to Bochner-Hartogs theorem \cite{Hor}, to $f$ being  the boundary values of a holomorphic function in $D.$

Our results on estimating CR dimensions of real manifolds, mentioned in the beginning of this subsection, lead to solutions of the above problem in real-analytic case.

Now let us explain the link between the result about CR dimension ( which in turn is a geometric form of Theorem \ref{T:Main})  and One-dimensional Holomorphic Extension Problem, in particular, with it special case-the strip-problem.

That link is almost immediate if one
passes from the function $f,$ for which one would like to test CR-condition, to its graph $\Lambda=graph_{\Omega} f.$
Then the analytic discs attached to $\Lambda$ are just the graphs of the analytic extensions of $f$
and CR-conditions for the function $f$ can be expressed in terms of CR-dimension of $\Lambda.$
To  apply the above formulated Parametric Argument Principle, we have to check Conditions 1 and 2. In the case of strip-problem Condition 1 holds because $H_2=0$ for any domain in the plane. Condition 2  turns  to the condition that the curves $\gamma_t$ surround no common point.

In higher dimensions, in the Globevnik-Stout conjecture, Condition 1 is satisfied automatically because the dimension of the manifold $\Lambda$ is less than $k+1,$ where $k$ is the dimension of the manifold of parameters.  Condition 2 holds because the analytic discs fill a domain with a hole $D^{\prime}.$

\bigskip

\noindent
{\bf Border problem}

The border problem is to determine whether a given real smooth closed manifold is the boundary of a complex manifold. This problem was solved in the works of Wermer and Harvey-Lawson
\cite{AW1, W1, W2, HL}.  It was proved there that necessary and sufficient condition for a real manifold $\Lambda \subset \mathbb C^n$  of dimension $2p-1$ to bound a holomorphic manifold is $dim_{CR}\Lambda=p$ when $p>1,$ and the
moment condition $\int_{\Lambda}\omega=0$ (for any holomorphic 1-form $\omega$), when $p=1.$

Furthermore, it is proved in articles of Dolbeault and Henkin \cite{DH0, DH} and of Dinh \cite{D1,D2,D3}, that the property of real $(2p-1)-$ manifold to bound a holomorphic manifold of dimension $p$ or, more generally,
holomorphic $p-$chain, can be verified slice by slice, i.e. via (one-dimensional) sections of the manifold by complex $(n+p-1)-$ planes $L \subset \mathbb C^n$. In this article, we prove a generalization of this result, by replacing linear sections to generic families of attached analytic discs.

\section{ Content of the article}

Section \ref{S:Basic} contains basic notions and notations we use throughout the article.

In Section \ref{S:parametric} we formulate Theorem \ref{T:Main} about parametric argument principle.

In Section \ref{S:Proof}  we prove Theorem \ref{T:Main}, which is the key result of the article.
Section \ref{S:Examples} contains examples.

Sections \ref{S:Lower_bounds} and \ref{S:Lower_general} are devoted to obtaining estimates from below of  CR-dimensions in terms of attached analytic discs, or equivalently, of complex moment conditions.
Theorem \ref{T:Thm2} is a geometric version of Theorem \ref{T:Main} and says that  if $d \leq n,$ then a real $d$-submanifold of a complex space has everywhere positive CR-dimension
(i.e., is nowhere totally real) if it admits a degenerate family of attached analytic discs with homologically nontrivial orbit.
Corollary \ref{C:Morera_curves} gives a characterization of complex curves in $\mathbb C^2$ as  real 2-manifolds admitting nontrivial families of attached analytic discs.

When the dimension $d$ of the manifold exceeds the dimension of the ambient space, i.e. $d>n$, then the CR dimension, which is  always at least $d-n,$ is positive and hence Theorem \ref{T:Thm2} becomes vacuous.
A nontrivial extension of this theorem to the range $d \geq n$  is given in Theorem \ref{T:Thm3}, where the  condition $\dim_{CR}\Lambda >0$ is derived in terms of vanishing of complex moments on curves, or, equivalently,
in terms of attached analytic discs.
Namely, Theorem \ref{T:Thm4} establishes the general  lower bound $\dim_{CR} \Lambda \geq q$ in terms of families of attached analytic discs. This estimate becomes
nontrivial  when $q > d-n.$

Theorem \ref{T:Thm5} and Theorem \ref{T:Thm6} are corollaries of \ref{T:Thm4}. They characterize, in terms of moment conditions,  manifolds or chains which are
holomorphic manifolds or boundaries of holomorphic manifolds, depending on the parity of the dimension. Such manifolds correspond to the cases of maximal possible CR-dimensions.

Section \ref{S:applications_Morera} is devoted to characterization of CR-functions defined on CR-manifolds. The obtained results are special
cases of corresponding theorems from Sections \ref{S:Lower_bounds} and \ref{S:Lower_general},  when the manifolds under consideration are graphs of functions.
We give characterizations CR-functions and boundary values of holomorphic functions in terms of analytic extension into families of
attached analytic disc (one-dimensional extension property).
In particular, these results  contain solutions of the strip-problem and of Globevnik-Stout conjecture for real-analytic functions, given earlier in
\cite{A1},\cite {A2}).

\section{Basic definitions and notations}\label{S:Basic}

We will use upper index to indicate the dimensions of manifolds. For real manifolds it will mean the real dimension, and
and for complex manifolds - the complex one. However, sometimes we will be omitting that index.

We start with some notions concerning CR manifolds. For more information about CR manifolds and functions we refer to a recent monograph \cite{MP} on the topic. 

Let $\Lambda=\Lambda^d \subset \mathbb C^n$ be a real smooth  $d-$dimensional submanifold, possibly with nonempty boundary.
For each  point $b \in \Lambda,$ we denote $T_b \Lambda$ the real tangent space to $\Lambda$ at the point $b$ and by
$$T_b^{\mathbb C}\Lambda= T_b(\Lambda) \cap i T_b^{\mathbb C}(\Lambda)$$-
the maximal complex tangent subspace.

Denote $$c(b)=c_{\Lambda}(b)= \dim_{\mathbb C} T_b^{\mathbb C}(\Lambda)$$
the complex dimension of the complex tangent space. The number $c(b)$ is called {\it the CR-dimension} of $\Lambda$ at the point $b.$ It is easy to see that
$$c(b) \ge d-n.$$

For $d \geq n$, the dimension $c(b)$ takes the minimal possible value $c(b)=d-n$ if and only if $T_b\Lambda$ generates a complex linear space of maximally possible dimension :
$$T_b\Lambda + i T_b\Lambda= T_b X.$$
In this case the manifold $\Lambda$ is called {\it generic.} When $d=n$ , the generic manifolds are {\it totally real}, i.e.,
have no nontrivial complex tangent spaces.
For $ d >n$ the manifold $\Lambda^d$ is never totally real, while when $d \leq n$ the manifold may be totally real or not.

When $c(b)=const$ then $\Lambda$ is called CR manifold. In this case, a smooth function $f$ on $\Lambda$ is called CR function
if it satisfies the tangential CR-equation,
$$\overline Zf=0$$
for any  complex tangent vector field $Z \in T^{\mathbb C}\Lambda.$

We  will be dealing with families of smooth analytic discs, $D_t$, smoothly parametrized by points $t$ belonging to a connected real-analytic oriented closed
$k-$ dimensional manifold, $M=M^k.$ The manifold $M^k$ can be assumed embedded in a Euclidean space $\mathbb R^N.$

Each analytic disc $D_t, \ t \in M,$ is parametrized by a holomorphic immersion
$$\Phi_t:\Delta \to D_t,$$
which is smooth up to the boundary of the unit disc.
The analytic disc $D_t$ is {\it attached} to a manifold $\Lambda \subset \mathbb C^n$ if $\partial D_t=\Phi_t(\partial \Delta) \subset \Lambda.$

Smooth  dependence  of the family $D_t$
on the parameter $t$ means that the mapping
$$\Phi : \overline \Delta \times M \to X,$$
defined by
$$\Phi(\zeta,t)=\Phi_t(\zeta)$$
belongs to $C^{r}(\overline \Delta \times M).$
The mapping $\Phi$ will be called  {\it parametrization} of the family $D_t$.

The order $r$ of smoothness can be finite, or $r=\infty,$ or $r=\omega$. In the latter case we deal with real-analytic family.
Throughout the paper the families of attached analytic disc will be real-analytic.

Let $\Lambda^d \subset \mathbb C^n$ be a real submanifold as above. The fact that the discs
$$\{D_t\}_{t \in M}$$ are attached to $\Lambda$ simply means that
$\Phi(\partial \Delta \times M) \subset \Lambda$. Moreover, we assume  that $\Lambda$ is covered by the closed curves $\partial D_t,$ i.e.
$$\Phi(\partial \Delta \times M)=\Lambda.$$
In this case we say that $\Lambda$ {\it admits} the family
$$ \mathcal F_M= \{D_t\}_{t \in M}.$$

The parametrization $\Phi$ will be assumed regular:
\begin{Definition}\label{D:regular} We will call the mapping $\Phi$ {\bf regular} if
\begin{enumerate}
\item
$rank \ d_t\Phi(\zeta,t) =\dim M, \ (\zeta,t) \in \overline \Delta \times M,$
where $d_t\Phi$ is the partial differential with respect to  $t.$
\item The image $\Lambda:=\Phi(\partial \Delta \times M)$ is a manifold, maybe with boundary, and
$rank \ d\Phi\vert_{\partial \Delta \times M}(u)$
equals $d=\dim \Lambda$ when $\Phi(u) \in \Lambda \setminus \partial \Lambda$ and  equals $d-1$ when $u \in \Phi^{-1}(\partial \Lambda).$
\end{enumerate}
\end{Definition}

Correspondingly, we will call the family $\mathcal F_M$ of the analytic disc {\it regular} if it admits a regular parametrization $\Phi.$


We will use also the following notations:
\begin{equation}\label{E:notation}
\Sigma =\Delta \times M, \  b\Sigma=S^1  \times M.
\end{equation}

\bigskip

\noindent
{\bf  I. MAIN RESULT}
\section{Formulation of Main Theorem} \label{S:parametric}

Most of results of this article are applications of the  following theorem about propagation of boundary degeneracy. We regard that theorem
as Parametric Argument Principle (see Section \ref{S:Motivation}).

Everywhere in the sequel, the homology groups are understood with coefficients in $\mathbb R.$

If $\Phi(\zeta,t)$ is a differentiable mapping from $\overline \Delta \times M$ to $\mathbb C^n,$ where $M$ is a differentiable $k$-dimensional  manifold and $\Phi(\zeta,t)$ is holomorphic in $\zeta \in \Delta$
then
$$rank_{\mathbb C}d\Phi(\zeta,t)=rank_{\mathbb C}\{\frac{\partial \Phi}{\partial \zeta}(\zeta,t),\frac{\partial\Phi}{\partial t_1}(\zeta,t),\cdots,\frac{\partial\Phi}{\partial t_k}(\zeta,t)\}$$
is the maximal number of linear independent , over $\mathbb C$, vectors in the system. Here $t_j$ are local coordinates on $M.$  Correspondingly,
$$rank_{\mathbb R}d\Phi(\zeta,t)=rank_{\mathbb R}\{\frac{\partial\Phi}{\partial x}(\zeta,t),\frac{\partial\Phi}{\partial_y}(\zeta,t), \frac{\partial \Phi}{\partial t_1}(\zeta,t),\cdots,\frac{\partial\Phi}{\partial t_k}(\zeta,t)\}$$
is the maximal number of $\mathbb R$-linear independent vectors in the system. Here $\zeta=x+iy.$

\begin{Theorem}(Parametric Argument Principle)\label{T:Main} Let $M=M^k$ be a  compact oriented real-analytic connected $k-$ dimensional closed manifold with the trivial homology group $H_{k-2}(M^k)=0.$
Let $\Phi:\overline \Delta \times M^k \to \mathbb C^n$ be a regular (see Definition \ref{D:regular} ) real-analytic mapping , holomorphic in the variable
$\zeta \in \overline \Delta,$ which maps  $b\Sigma=S^1 \times M^k$ onto a real-analytic variety $\Lambda=\Phi(b\Sigma)$ of dimension $ d \leq n.$
Suppose also that
\begin{enumerate}
\item
The induced homomorphism of
the  homology groups
$$\Phi_*:H_{k+1}(b\Sigma ) \to H_{k+1}(\Lambda)$$
is zero, $\Phi_*=0.$
\item
The induced homomorphism of the  homology groups
$$\Phi_* : H_k(\overline \Sigma) \to
H_k(\Phi(\overline \Sigma))$$
is not zero, $\Phi_* \neq 0.$
\end{enumerate}
Then $rank_{\mathbb C}(\zeta,t) < d$ for all  $(\zeta,t) \in \overline \Delta \times M^k.$
None of the  conditions 1 and 2 can be omitted.
\end{Theorem}

\begin{Remark}

\begin{itemize}
\item
The condition $H_{k-2}(M^k)=0$ required by Theorem \ref{T:Main} is a technical one, needed
in the proof of Theorem \ref{T:Main}.
\item
The conclusion of Theorem \ref{T:Main} means that $\Phi$ degenerates in the sense of its rank. Indeed, since $\Phi(S^1 \times M^k)=\Lambda,$ we have $d=dim \Lambda \leq k+1.$ On the other hand, the complex rank of the system
\begin{equation}\label{E:deriv}
\{\partial_x \Phi, \partial_{y} \Phi, \partial_{t_1}\Phi,\cdots,\partial_{t_k}\Phi\}
\end{equation}
of $k+2$ vectors  is at most $k+1,$ because $\partial_{x}\Phi=i\partial_{y}\Phi$ due to Cauchy-Riemann equations.  Therefore, the inequality $rank_{\mathbb C} d \ \Phi < d \leq k+1$ implies that the rank of $\Phi$ is everywhere less than the maximal possible one, which is $k+1.$
It occurs, in particular, when $\Phi$ maps the $k+2$-dimensional manifold $\overline\Delta \times M^k$ into a manifold of a less dimensions, i.e., $rank_{\mathbb R} d \Phi(\zeta,t) <k+2, \ (\zeta,t) \in \overline \Delta \times M^k.$ Indeed,  then the system (\ref{E:deriv}) is $\mathbb R-$ linearly dependent, and the rank of that system over $\mathbb C$ is less than $k+1,$  due to CR equation $\partial_x \Phi=i\partial_y \Phi.$
\end{itemize}
\end{Remark}

Notice also, that since, on one hand, $d=\dim \Lambda=\dim \Phi(b\Sigma) \geq k+1$ and, on the other hand,  $d \geq k,$ by regularity condition, the dimension $d$ can take two values: either $d=k$ or $d=k+1.$
\begin{Proposition}\label{P:condition_1}

Condition 1 in Theorem \ref{T:Main} can be reworded as follows:

\noindent
a) either $\Phi$ degenerates on $b\Sigma$ in the sense of dimension, i.e.
$d < k+1$

\noindent
or

\noindent
b) $d=k+1$
and $\Phi$ topologically degenerates on
$b\Sigma,$ meaning that  its topological (Brouwer) degree (cf. \cite{Hirsch}, \cite{Milnor}) is zero: $deg \ \Phi\vert_{b\Sigma}=0.$

\end{Proposition}
\pf
The $(k+1)$-st homology group of the $(k+1)$-manifold $b\Sigma=S^1 \times M^k$ is generated by the fundamental homology class
$$[b\Sigma]=[S^1 \times M^k].$$
The maximal rank of $\Phi$ on $b\Sigma$ is $k+1,$ while  the dimension of the target manifold $\Lambda$ is at most $k+1.$
If $\Lambda$ has dimension $k$ or  less, then $\dim \Lambda <k+1$ and hence the $(k+1)$-th homologies of $\Lambda$ are trivial and
$\Phi_*=0$. Otherwise,
$\Lambda$ and $b\Sigma$ have the same dimension $k+1$ and then
$$\Phi_*[b\Sigma]=\deg \Phi \cdot [\Lambda],$$
where $\deg \Phi$ is the Brouwer
degree. Since due to condition  $\Phi_*[b\Sigma]=0,$ we have $\deg \Phi=0.$
\qed

\begin{Remark}\label{R:k=0}
If we take $k=0$ in Theorem \ref{T:Main}
then $M$ degenerates to a point,  $\Sigma=\Delta$, and $\Phi$ is an analytic function  in the unit complex disc $\overline \Delta.$ Also, $\Lambda=\Phi(\partial \Delta)$ and  the dimension $d=dim \Lambda$ is either $d=0$ or $d=1.$ Condition 1
of Theorem \ref{T:Main} in this case means that either $\Lambda$ is a point, or $\Lambda$ is a curve with nonempty boundary so that the mapping $\Phi:S^1 \to \Lambda$
has the topological degree 0. However, the latter possibility is excluded.
Indeed, the conclusion of our theorems  says that $rank_{\mathbb C} d \Phi < d \leq 1$ and therefore $rank_{\mathbb C} d\Phi=0$
which means that $\Phi=const,$ in accordance with the classical argument principle.

Thus, Theorem \ref{T:Main} can be viewed as a generalization of the argument principle from a single analytic function $\Phi(\zeta)$
in $\Delta$ to a family $\Phi_t(\zeta)=\Phi(\zeta,t)$ of analytic functions in the unit discs, where parameter $t$ runs a real
manifold $M^k.$
\end{Remark}
We conclude this section with agreeing on the terminology for Conditions 1 and 2, which we will be using in the sequel.
\begin{Definition}\label{D:degen} We say that the mapping $\Phi:\Delta \times M^k \to \mathbb C^n$ {\bf (homologically) degenerates} on the boundary if
Condition 1 from Theorem \ref{T:Main} is fullfiled.
\end{Definition}
\begin{Definition} \label{D:nontrivial} We say that the mapping $\Phi:\Delta \times M^k \to \mathbb C^n$ has {\bf homologically nontrivial $t$-orbit}
if Condition 2 from Theorem \ref{T:Main} is fullfiled. This means that the image $\Phi(\{\zeta_0\} \times M^k)$
of the fundamental cycle (the orbit), for some (and then for any) point $\zeta_0 \in \overline\Delta$   is
not homological to 0 in $\Phi(\overline\Delta \times M).$ The latter means that $\Phi(\{\zeta_0\}\times M)$ bounds no $k+1-$ chain in $\Phi(\overline\Delta\times M).$
\end{Definition}

\section{Proof of Theorem \ref{T:Main}}\label{S:Proof}

In this section we prove our main result, Theorem \ref{T:Main}, from which we  will derive all other results presented in Sections 7-11.
The proof develops the main idea from   \cite{A1}, \cite{A2}.

\subsection{Model of the proof (the case $n=1, k=0.$)}\label{S:Model}
To warm up, we will first demonstrate the idea of the proof in the simplest case when $n=1$ and the parameter $t$ is absent, that is when one deals with just a holomorphic function.
In this case, Theorem \ref{T:Main} reduces to Proposition \ref{P:argument_principle}.

We want to present an alternative proof of Proposition \ref{P:argument_principle}, which
is perhaps is not as natural as one presented in Section \ref{S:Motivation} and even  requires analyticity in the closed unit disc, but instead contains a hint to generalization to the parametric case.
\begin{Proposition}\label{P:arg_principle_real_a} Let $\Phi$ be a function holomorphic in the closed unit disc $\overline\Delta$ and $\gamma:=\Phi(\partial\Delta).$ If the induced homomorphism $\Phi:H_1(\partial\Delta) \to H_1(\gamma)$ is zero, then $\Phi^{\prime}=0$, i.e., $\Phi=const.$
\end{Proposition}
\pf Denote $J(\zeta)=\partial_{\psi}\Phi(re^{i\psi})=i\zeta \Phi^{\prime}(\zeta), \ \zeta=re^{i\psi}.$ We need to prove that $J=0$ identically.

Suppose that this is not the case.
Then $J$ has finite number of isolated zeros in $\overline\Delta.$ If $\Phi$ maps two points $\zeta_1,\zeta_2 \in S^1=\partial\Delta$ to the same point $b:=\Phi(\zeta_1)=\Phi(\zeta_2)$ in $\gamma$ and $J(\zeta_1), J(\zeta_2) \neq 0$, then
the vectors $\partial_{\psi} \Phi(\zeta_1)$ and $\partial_{\psi} \Phi(\zeta_2)$ are tangent to $\gamma$ at the point $b$ and hence they are collinear, $\partial_{\psi} \Phi(\zeta_1)=\lambda \partial_{\psi} \Phi(\zeta_1), \lambda \in \mathbb R.$
Then we have
$$\frac{J\zeta_1)}{\overline {J(\zeta_1)}}=\frac {J(\zeta_2)}{\overline {J(\zeta_2)}}.$$
The identity holds for all points $\zeta_1,\zeta_2 \in \partial \Delta$ except for finite number of zeros of $J.$

Due to real-analyticity, then $J/\overline J$ represents on $\partial\Delta \setminus J^{-1}(0)$ as
\begin{equation}\label{E:sig}
\frac{J}{\overline J}=\sigma \circ \Phi,
\end{equation}
where $\sigma$ is a smooth (real-analytic) function on $\Phi(\partial \Delta \setminus J^{-1}(0)).$
Now consider the logarithmic residue integral:
\begin{equation}\label{E:log}
I=\frac{1}{2 \pi i}\int\limits_{\partial \Delta}\frac{d(J/\overline J)}{J/\overline J}= \frac{1}{2 \pi i}( \int\limits_{\partial \Delta}\frac{dJ}{J}- \int\limits_{\partial \Delta}\frac{d\overline J}{\overline J}).
\end{equation}
Due to real-analyticity of $J$, the integrals converge at zeros of $J$ in the sense of principal value.

The function $J(\zeta)$ is holomorphic in $\overline \Delta$ and hence,  according to Argument Principle, the right hand side in (\ref{E:log}) computes the number of zeros of $J$ in the closed disc, counting multiplicities
(the boundary zeros contribute with the factor $\frac{1}{2}).$  Since $J$ has zero at the origin $\zeta=0,$ we have  $I>0.$
On the other hand, using (\ref{E:sig})  we have for the left hand side in (\ref{E:log}) after change of the variables $z=\Phi(\zeta)$:
$$I=\int\limits_{\partial\Delta} \frac{d(\sigma \circ \Phi)}{\sigma \circ \Phi}= deg \ \Phi\vert_{\partial\Delta}\int\limits_{\gamma}\frac{d\sigma}{\sigma}=0$$
because $\deg \ \Phi=0$ on $\partial \Delta$ by the condition. This contradiction proves that $J=0$ identically and hence $\Phi^{\prime}=0.$
\qed

The proof of Theorem \ref{T:Main} follows the strategy of the above proof of  Proposition \ref{P:arg_principle_real_a}.

First, we reduce the conclusion of Theorem \ref{T:Main} to proving that certain Jacobians (analogs of the derivative $J$ in the above example) vanish identically.

Assuming that this is not the case, we construct a nontrivial cycle (orbit of the origin $\zeta_0=0$) in $\Delta \times M^k$ consisting of critical points of $\Phi$ -- zeros of $J.$ By the condition, $\Phi$ maps this cycle to a nontrivial cycle in the image.

By the de Rham duality, there exists a closed differential form  having nonzero integral over the image cycle.

Applying  Stokes formula, which relates that integral with the integral over
the boundary $\partial\Delta \times M^k$, we arrive at contradiction, because the integral over the boundary is zero due to Condition 1 of homological degeneracy of $\Phi$ on the boundary.

\subsection{Preparations}

Let the mapping $\Phi:\overline \Sigma=\overline \Delta \times M^k \to \mathbb C^n$
be as in Theorem \ref{T:Main} and suppose that all the conditions of Theorem \ref{T:Main} are fullfiled.

\begin{Lemma} \label {L:L1}
Fix two natural numbers $1< d \leq k+1.$ Suppose that for any holomorphic $d$-form $\eta$ in $\mathbb C^n$ and for any $d-1$
smooth (real-analytic) tangential vector fields $T_1,\cdots,T_{d-1}$ on the manifold $M^k$ the following identity holds:
\begin{equation}\label{E:Jacobi}
J(u):=\eta(u; \partial_{\psi}\Phi(u), (T_1\Phi)(u),\cdots,(T_{d-1}\Phi)(u))=0, \ u \in b\Sigma:=S^1 \times M^k.
\end{equation}
Then  $rank_{\mathbb C} \ d\Phi(u) < d , \ u \in \overline \Sigma.$
\end{Lemma}

Here $T_j$ are globally defined tangential vector fields on $M^k$, which are allowed to vanish on less dimensional sets in $M,$ and which are
viewed as tangential differential operators
$$T_j=\sum\limits_{i=1}^k a_{ij}(\varphi)\frac{\partial}{\partial\varphi_i}$$
on $M^k,$ $\varphi_i$ are local coordinates on $M^k$: $t=t(\varphi_1,\cdots,\varphi_k).$

\pf First of all, observe that there are sufficiently many real-analytic tangential
vector fields on the real-analytic manifold $M,$ in the following sense: for every point
$t_0 \in M^k$ there exist $k$ real-analytic tangential vector fields $T_1,...,T_k$ on
$M^k$ such that $T_1(t_0),...,T_k(t_0)$ form a basis in the tangent space $T_t(M^k).$

Indeed, we can think of $M^k$ as a submanifold embedded to a Euclidean space $\mathbb
R^N$ and take any orthonormal basis $e_1(t_0),...,e_k(t_0)$ in $T_{t_0}(M^k).$ If
$$P_t:\mathbb R^N \to T_t(M^k):$$
is the operator of orthogonal projection, then for every $x \in \mathbb R^N$ the vector
field $P_t(x),  \ t \in M^k$ is real-analytic, since near each point it expresses as
$$P_t(x)=\sum_{j=1}^k <x,e_j(t)>e_j(t),$$
where $e_j(t), \ j=1,...,k$ is an orthonormal tangential frame locally defined in a
parametrized neighborhood and real-analytically depending on $s.$

Then the tangential vector fields $P_t(e_1(t_0)),...,P_t(e_k(t_0))$ are what we need.

Now consider the two possible cases:

\noindent
1) $d-1=k.$
The condition of Lemma says that for any $u=(\zeta,t)=(e^{i\psi},t) \in b\Sigma$  the vectors
$$\partial_{\psi}\Phi(u), T_1\Phi(u),\cdots,T_k\Phi(u),$$
are linearly dependent over $\mathbb C.$ Since, as it was shown, for any point $t \in M$
the vector fields $T_j$ can be chosen so that they form a basis at $T_t(M),$  the
condition exactly means that $rank_{\mathbb C} \ d\Phi(\zeta,t) <k+1=d$ for $|\zeta|=1.$
Since $\Phi(\zeta,t)$ is holomorphic with respect to $\zeta$, the inequality extends for
$|\zeta|<1.$

\noindent 2) $d \leq k.$ Take  $u=(e^{i\psi},t) \in b\Sigma$ and let $T_j, \ j=1,...,k$
constitute a basis in $T_t(M^k).$
 If $rank_{\mathbb
C} \ d\Phi(u) \geq d$ then the inequality holds for $u$ in  an open set $U \subset b
\Sigma.$ This means that the system
$$\partial_{\zeta}\Phi(u)=-ie^{-i\psi}\partial_{\psi}\Phi(u), T_1\Phi(u),\cdots,T_k \Phi(u), \ u \in \overline \Sigma$$
contains at least $d$ linearly independent, over $\mathbb C,$ vectors.

By condition ( \ref{E:Jacobi}) $\partial_{\psi}\Phi(u)$ is not one of these $d$ vectors.
Therefore, the subsystem $\mathcal E$ of $d$ linearly independent vectors consists of derivatives with respect to the parameters $t_j$ only:
$$\mathcal E=\{T_{i_1}\Phi(u),\cdots,T_{i_d}\Phi(u)\}.$$

Any $d-1$ vectors in $\mathcal E$ are also linearly independent. Condition (\ref{E:Jacobi}) says that adding $\partial_{\psi}\Phi(u)$  to any system of $d-1$ elements from $\mathcal E $ makes the system $\mathbb C-$ linearly dependent.  Therefore,
$\partial_{\psi}\Phi(u)$ belongs to  any $\mathbb C$-linear subspace spanned by  $d-1$ vectors in $\mathcal E.$ The intersections of those linear subspaces is $0$ and therefore $\partial_{\psi}\Phi(u)=0.$

Since $u$ belongs to an open set $U$ in $b\Sigma=\partial \Delta \times M^k$ then by real-analyticity $\partial_{\psi}\Phi(\zeta,t)=0$ for $|\zeta|=1.$ Since $\Phi(\zeta,t)$ is holomorphic in $\zeta,$ we conclude that $\Phi$ does not depend on $\zeta,$ which is not the case.

Thus, $rank_{\mathbb C}d\Phi(u) <d$ for all $u \in \partial \Delta \times M^k.$ This means that all $d\times d$ minors of the Jacobian matrix of $J(\zeta,t)$ with respect to $\zeta,T_1,\cdots,T_k$ vanish for $|\zeta|=1$.
Since $J$ is holomorphic in $\zeta$, it is so for $|\zeta|<1$ and hence $rank_{\mathbb C}d\Phi(u)<d$ for all $u \in \overline\Delta\times M^k.$ Lemma is proved.

\pf

\subsection{The Jacobians}\label{S:Jacobians}

Suppose that the conclusion of Theorem \ref{T:Main}, which we are going to prove, is not true, i.e.,  the mapping $\Phi$ is nondegenerate in the sense of its rank.

Then  Lemma \ref{L:L1} yields that
there exists a holomorphic $d$-form $\eta$ and  $d-1$ smooth tangential vector fields $T_1,\cdots,T_{d-1}$ on $M^d$ such that the function
\begin{equation}\label{E:definitionJ}
J(u)=\eta(u;\partial_{\psi}\Phi(u),T_1\Phi(u),\cdots,T_{d-1}\Phi(u))
\end{equation}
does not vanish identically  for $u=(\zeta,t) \in S^1 \times M^n.$ Our goal is to arrive to contradiction.

From now on to the end of the proof, we fix the above form $\eta$ and the vector fields $T_j, \ j=1,\cdots,k.$
Observe that the function $J$ is naturally defined in the solid manifold $\overline \Delta \times M^n$
because $\Phi(\zeta,t)$ and its derivative in the angular variable, $\partial_{\psi}=i\zeta\partial_{\zeta}$,
are defined for $|\zeta| \leq 1.$
The differential operators $T_j$ act only in the variable $t \in M^n$ and do not depend on $\zeta$.

We will call the function $J$ {\it Jacobian}. Our main assumption is that $J$ is not identically zero and our aim is to obtain a contradiction
with this assumption.

\begin{Lemma} \label{L:Jacob}
The Jacobian $J(\zeta,t) $ has the following properties:

\noindent
a) $J(\zeta,t)$ is holomorphic in $\zeta \in \overline \Delta.$

\noindent
b) $J(0,t)=0, \forall t \in M^n.$

\noindent
c) the function $J/\overline J$ can be represented on $(S^1 \times M^n) \setminus J^{-1}(0)$ as
\begin{equation}\label{E:Jacob}
\frac{J}{\overline J}=\sigma \circ \Phi,
\end{equation}
for some smooth (real-analytic)  function $\sigma$ on $\Lambda=\Phi((S^1 \times M^k) \setminus J^{-1}(0)).$

\end{Lemma}
\pf
The property a) follows from the definition (\ref{E:definitionJ}) of the Jacobian because $\Phi(\zeta,t)$ is holomorphic in $\zeta$ and
$\eta$ is a holomorphic form. The property b) is due to vanishing  the derivative
$$\partial_{\psi}\Phi(\zeta,t)=i\zeta\partial_{\zeta}\Phi(\zeta,t)$$
at $\zeta=0.$

Let us prove the property c).
Take any two points $u_1, u_2 \in S^1 \times M^n$ such that
$\Phi(u_1)=\Phi(u_2).$

The condition $J(u_1), J(u_2) \neq 0$  implies that each of the two systems of $d$ vectors
$$\mathcal E^{j}=\{\partial_{\psi}\Phi(u_{j}), T_1\Phi(u_j,\cdots, T_{d-1}\Phi(u_j)\} , \ \ j=1,2,$$
is $\mathbb R-$ linearly independent.  Since all these vectors belong to the same $d$-dimensional space $T_b(\Lambda)$,
each of the two systems  constitutes a basis in this space.
The transition $d \times d-$ matrix $A,$ from the basis $\mathcal E^1$ to the basis $\mathcal E^2,$ is real and the Jacobians
differ by the determinant of the transition matrix:
$$J(u_1)=\det A \cdot J(u_2).$$
Since the determinant of $A$ is  real- valued , we have
$$\frac {J(u_1)}{\overline {J(u_1)}}= \frac {J(u_2)}{\overline {J(u_2)}},$$
 whenever $\Phi(u_1)=\Phi(u_2)$ and $J(u_1), J(u_2) \neq 0.$
For each $b \in \Lambda,$  which is  not a critical value, i.e., $b \notin \Phi(J^{-1}(0)),$ define
$$\sigma(b)=\frac{J(u)}{\overline J(u)}, \ \Phi(u)=b.$$  Then the property c) follows.
\qed

\subsection{The structure of the critical set $J^{-1}(0)$ and the current $d\ln J$. }\label{S:critical_set}

Let us examine the structure of the critical set $J^{-1}(0).$ We want to show that this set
determines $k$-current in $\overline \Delta \times M^k$ by means of the singular form $$d\ln J=\frac{dJ}{J}.$$
The value of this $k$-current at a $k$-form $\omega$ is defined by
integration of the wedge-product
$$d\ln J \wedge \omega.$$  This operation, in turn,  evaluates the integral of $\omega$ over the zero set of $J,$ multiplied by $2\pi i.$
The justification of this relation is given in Appendix \ref{S:Appendix}(see also \cite{A1}).

Since $J$ is real-analytic, the equation $J(\zeta,t)=0, \ (\zeta,t) \in \Sigma=\overline \Delta \times M^k$
defines on $\Sigma$ an analytic set \cite{F},\cite{GM}.
The function $J(\zeta,t)$ is analytic with respect to $\zeta$ in the closed unit disc, hence for each fixed $t \in M^k$ the
function $J(\zeta,t)$ either has finite number of isolated zeros in $\overline \Delta$ or vanishes  identically.
Correspondingly,  the zero set of $J$ can be decomposed into two parts:
$$J^{-1}(0)= N_{reg} \cup N_{sing},$$
where $N_{reg}$ is defined by the isolated zeros $\zeta=\zeta_j(t) \in \Delta$ and is $k-$ dimensional, and $N_{sing}$ consists of those discs
where $J$ vanishes identically with respect to $\zeta:$
$$N_{sing}=\{(\zeta,t):J(\zeta,t)=0, \forall \zeta \in \overline\Delta\}.$$

The regular $k$-dimensional part $N_{reg}$ can be regarded as the $k$-dimensional chain,
which is the sum of the singular chains defined by the mappings $$t \to
(\zeta_j(t),t), \ t \in U \},$$ where $t$ runs a simplex $U$ from a simplicial partition of $M^k.$ Here
$\zeta_j(t)$ is a local branch of (isolated) zeros of the holomorphic function $\zeta \to
J(\zeta,t).$

The number of zeros, counting multiplicities, near each inner zero $\zeta_j(t_0)$ is
computed by the logarithmic residue
$$\frac{1}{2\pi}\int\limits_{|\zeta-\zeta_j(t_0)|=\delta}\frac{dJ(\zeta,t)}{J(\zeta,t)}$$
which  continuous depends on $t$ so along as $J$ does not vanish on the circle $|\zeta-\zeta_j(t)|=\delta.$ Therefore the number of zeros is constant when $\delta$ is sufficiently small and $t$ is sufficiently close
to $t_0.$ 

The multiplicty assigned to a singular chain $t \to (\zeta_j(t),t), \ t \in U$  equals to the multiplicity of the
zero $\zeta_j(t)$ which is constant on each singular chain. The brancing points, at which  some zeros $\zeta(t)$ merge, constitute  a less dimensional subvariety of $M.$  The same logarithmic reside integral shows that at branching point the multiplicities of merging chains add.

Also, the boundaries of the branches lie on the unit circle $|\zeta|=1.$
The orientation is defined by the orientation of the parameter space $M^k.$ The
integration of $k-$ forms over $N_{reg}$ is therefore well defined.

 The second, singular, part $N_{sing}$ decomposes to the direct product:
$$N_{sing}=\overline \Delta \times T,$$
where $T \subset M^k$ is an analytic subset.

In turn, the analytic set $T$ falls apart into strata $T^s$  of different dimensions (see \cite{GM}).
Correspondingly, the singular set $N_{sing}$ falls apart into $s+2-$ dimensional strata $\overline \Delta \times T^s:$
$$N_{sing}=\cup_s (\overline \Delta \times T^s).$$
The strata with $s \leq k-3,$ are negligible from the point of view of integration of $k$-forms,  because
they have the dimension $s+2 < k.$  So, only the strata $\overline\Delta \times T^s$ with $s=k-2$ and $s=k-1$  can contribute to the integrals of $k-$ forms.

Furthermore, the singularities of the form $\overline \Delta \times T^s$ with $s=k-1$ are
removable  due to the following argument. The stratum $T^{k-1}$ is locally defined by
zeros of a real function, $T^{k-1}=\{\rho=0\},$ and then  $J(\zeta,t)=\rho(t)^m \
J_0(\zeta,t)$  for some positive integer $m,$ where $J_0(\zeta,t)$ does not vanish
identically with respect to $\zeta$ when $t \in T^{k-1}.$ Then the the real factor
$\rho^m$ cancels in the ratio
$$\frac{J}{\overline J}=\frac {I}{\overline J_0}.$$ This ratio is exactly what we are going to integrate in the sequel, so that the singularities corresponding to $T^{k-1}$ are removable with respect to integration of the function $\frac{J}{\overline J}.$


Thus, only the $k$-dimensional part
$$N_{reg}=\overline \Delta \times T^{k-2},$$
corresponding to the strata $T^{k-2}$ of dimension $s=k-2,$
can contribute in the current $d\ln J.$  However, the condition $H_{k-2}(M^k)=0$ of Theorem \ref{T:Main} excludes such a possibility as well.

The analytic justification of manipulations with the logarithmic derivative $d ln J,$ viewed as a current, is given in Appendix \ref{S:Appendix}.

\subsection {The duality}\label{S:duality}
We refer for the details of homological duality theory to the books \cite{S}, Ch.6, Section 2, and \cite{Hatcher}, Section 3.3.


{\bf The dual differential form $\omega.$} Let $Y$ be a $m-$ dimensional compact manifold. The cohomology group $\overline{H}^k(Y)$ is dual to the homology group $H_k(Y)$ and hence if the homology class
$[C] \in H_k(Y)$ and $[C] \neq 0$ then there exists a cocycle $h \in H^k(Y)$ such that the pairing $\langle [C], h \rangle \neq 0.$
According to the de Rham realization of the  cohomology group $H^k(Y, \mathbb R)$ by differential forms,   $h$ is represented by
 a closed differential $k-$ form $\omega$ such that the pairing
$$\langle [C], \omega \rangle =\int\limits_C \omega \neq 0.$$

\subsection {Main Lemma}

The proof of Theorem \ref{T:Main} exploits the above mentioned fact that the current defined by the differential form $d lnJ=\frac{dJ}{J}$ corresponds to the integration over the critical set $J^{-1}(0).$
The following lemma is the key point of the proof.

\begin{Lemma}\label{L:stokes}
For any closed differential $k$-form $\omega$ on $\Phi(\overline\Delta \times M^k)$
holds
\begin{equation}\label{E:7}
\int\limits_{J^{-1}(0)} \Phi^*\omega=0.
\end{equation}
In particular, it follows that the chain $\Phi(J^{-1}(0))$ is a cycle in $\Phi(\overline\Delta \times M^k).$
\end{Lemma}
\pf
First, check the last assertion. By change of variables, we have from (\ref{E:7})
$$\int\limits_{\Phi(J^{-1}(0))}\omega=0$$
for any closed $k-$ form $\omega.$ In particular, $\omega$ can be taken as $\omega=d\eta,$ where $\eta$ is an arbitrary $k-1$-form. Then by Stokes formula
$$\int\limits_{\Phi(J^{-1}(0))}\omega=\int\limits_{\partial \Phi(J^{-1}(0))}\eta=0.$$
Since  $\eta$ is an arbitrary $(k-1)-$ form, we conclude that $\partial\Phi(J^{-1}(0)) =\emptyset.$

Now,  $\omega$ is a closed form,  the differential form $d \ ln \frac{J}{\overline J} \wedge \Phi^*\omega $ is  closed outside of the critical set $J^{-1}(0)$ and hence
$$\int\limits_{(\overline \Delta \times M^k)  \setminus J^{-1}(0)} d(d\ln \frac{J}{\overline J} \wedge \Phi^*\omega)=0.$$
On the other hand, the Stokes formula for the currents yields
\begin{equation} \label{E:8}
\int\limits_{(\overline \Delta \times M^k)  \setminus J^{-1}(0)} d(d\ln \frac{J}{\overline J} \wedge \Phi^*\omega)=
\int\limits_{(S^1 \times M^k) \setminus J^{-1}(0) } d\ln \frac{J}{\overline J} \wedge \Phi^*\omega + 4\pi i \int\limits_{J^{-1}(0)} \Phi^*\omega.
\end{equation}
An analytic justification of this formula is given in Appendix \ref{S:Appendix}.

It remains to prove that the first summand
$$I:= \int\limits_{(S^1 \times M^k) \setminus J^{-1}(0)} d \ ln \frac{J}{\overline J} \wedge \Phi^*\omega$$
in the right hand side of (\ref{E:8}) is 0.

Denote
$$\Xi:=d\ln\frac{J}{\overline J} \wedge \Phi^*\omega.$$
Using the representation (\ref{E:Jacob})
from Lemma \ref{L:Jacob} for $J/\overline J $  we obtain
\begin{equation}\label{E:10}
\Xi=\Phi^* (d \ ln \sigma) \wedge \Phi^*\omega=\Phi^*(d \ln\sigma \wedge \omega),
\end{equation}
so that the integral $I$ can be written as
\begin{equation}\label{E:9}
I=\int\limits_{(S^1 \times M^k)\setminus J^{-1}(0)}\Xi=\int\limits_{(S^1 \times M^k)\setminus J^{-1}(0)} \Phi^*(d \ ln\sigma \wedge \omega).
\end{equation}

According to Proposition \ref{P:condition_1} the condition 1 of Theorem \ref{T:Main} implies that one of the  two cases is possible:

\noindent
a) $d <k+1,$

\noindent
b) $d=k+1$ and the Brouwer degree $\deg \ \Phi$ on $S^1 \times M^k$ equals 0.

Consider the case a), when  $k \geq d.$
Then the degree of the form $\omega$ is

$$ \deg \omega=k \geq d$$
and hence
$$\deg (d\ln\sigma \wedge \omega) =k+1 \geq d+1.$$

Since
$\dim \Phi(S^1 \times M^k)=\dim \Lambda^d=d<d+1,$  we have:
$$d\ln\sigma \wedge \omega\vert_{\Phi(S^1 \times M^d)}=0.$$
But  then the pull-back differential form $\Xi=\Phi^*( d\ln\sigma \wedge \omega)=0$ and hence
$$I=\int_{S^1 \times M^k} \Xi=0.$$

\medskip
Now turn to  the case b). In this case
$k=d-1,$ $\Phi$ maps $d-$ dimensional manifolds $S^1\times M^k$ to $d-$dimensional manifold $\Lambda^d$ and the Brouwer degree of this mapping $\deg \Phi=0.$
Then the integral (\ref{E:9}) vanishes.
Indeed,  the definition (\ref{E:10}) of $\Xi$ and change of variables yield
$$I=\int\limits_{(S^1 \times M^{d-1}) \setminus J^{-1}(0)} \Xi=
\deg \Phi \  \int\limits_{\Phi [(S^1 \times M^{d-1}) \setminus J^{-1}(0)]} d\ln \sigma \wedge \omega=0,$$
due to $\deg \Phi=0.$ Lemma is proved.
\qed

\subsection {End of the proof of Theorem \ref{T:Main}}

Now we are ready to finish the proof of Theorem \ref{T:Main}.  By Condition 2 , the image  $\Phi(\{0\} \times M^k)$ of the fundamental cycle $\{0\}\times M^k$ represents nonzero homological class
in the $k-th$ homology group of $\Phi(\overline \Delta \times M^k).$

By the de Rham duality (see subsection \ref {S:duality}) there exists a closed nonexact $k$-form $\omega$
such that
$$\int\limits_{\Phi(\{0\} \times M^k)} \omega \neq 0.$$
By changing the sign of $\omega$, we can assume that the integral is strictly positive.

As above, the cycle of integration consists of the two parts:
$$\Phi(J^{-1}(0))=\Phi(N_{reg}) \cup \Phi(N_{sing}).$$

The first part contains the image $\Phi(\{0\} \times M^k)$ of the fundamental cycle in
$\Delta \times M^k.$ The entire cycle $\Phi(N_{reg})$ is homological to $m \Phi(\{0\}
\times M^k)$ with $m> 0,$ because $\Phi(\zeta,t)$ is holomorphic in $\zeta$ and hence
preserves the orientation of the chains $(\zeta_j(t),t), t \in M^k$ constituting
$J^{-1}(0)$ (one of these chains is just the cycle $\Phi(\{0\} \times M^k)$.) See
\ref{S:Appendix} for more justification.

The second part $\Phi(N_{sing})$ is negligible, as it was shown in Section \ref{S:critical_set}.

Therefore, since $\omega$ is closed, we have
\begin{equation}\label{E:11}
\int\limits_{J^{-1}(0)}\Phi^*\omega=
\int\limits_{\Phi(J^{-1}(0))}\omega  \geq \int\limits_{\Phi(\{0\} \times M^k)}\omega >0,
\end{equation}

in  contradiction with  Lemma \ref{L:stokes}. The source of the contradiction is the assumption that $J$ is not identically zero. Thus, $J=0$ and and this completes the proof,
due to Lemma \ref{L:Jacob} and to the definition of the Jacobian $J$ in Section \ref{S:Jacobians}. Now the proof of Theorem \ref{T:Main} is complete.

\subsection{A justification of formula (\ref{E:8})}\label{S:Appendix}
As it was discussed in Section \ref{S:critical_set}, the only part of the critical set
$J^{-1}(0)$ that contributes in integration of differential $k-$ forms is the regular
part $N_{reg} \subset J^{-1}(0).$ 

This part is regarded, 
as it was explained in the beginning of Section \ref{S:critical_set},  as  a $k-$ chain which is the sum of
the singular chains defined by the mappings $t \to \zeta_j(t),$ where $\zeta_j(t)$ are finitely many
zeros of the holomorphic function $\zeta \to J(\zeta,t)$ in the closed unit disc, and $t$ runs a simplex in $M$ from a simplicial
partition. The singular chains participate in the integration of differential forms as the oriented manifolds parametrized
by $t \to (\zeta_j(t),t)$ with the multiplicty equal to the multiplicty of the zero $\zeta_j(t).$

Apply Stokes formula to the manifold
$$\Sigma_{\varepsilon}:=(\overline \Delta \times M^k) \setminus N_{reg,\varepsilon},$$
where $N_{reg,\varepsilon}$ is the tubular $\varepsilon-$ neighborhood of $N_{reg},$
defined for each singular chain $N_j$ in $N_{reg}$ as
$$N_{j,\varepsilon}:= \{(\zeta,t): |\zeta -\zeta_j(t)| <\varepsilon\},$$
and to the closed differential $k-$ form
\begin{equation}\label{E:J/J}
\frac {d(J/\overline J)}{J/\overline J} \wedge \Phi^*\omega=\frac{dJ}{J}\wedge \Phi^*\omega -
\frac{d\overline J}{\overline J}\wedge \Phi^*\omega,
\end{equation}
defined out of the critical set $J^{-1}(0).$
The orientation of $N_{reg,\varepsilon}$ is naturally defined by the orientations on $M^k$ and on the $\zeta-$ plane.

Since $d \ln\frac{J}{\overline J}=d\ln J -d\ln\overline J,$
the integral in (\ref{E:8}) splits into the difference of two integrals. We will apply Stokes formula separately to the integral containing $ d\ln J=dJ/J$ and one containing
$d \ln \overline J=d\overline J/\overline J.$

We have for the first integral containing $dJ/J:$
\begin{equation} \label{E:A1}
0=\int\limits_{\Sigma_{\varepsilon}} d(\frac{dJ}{J} \wedge \Phi^*\omega)=
\int\limits_{\partial \Sigma_{\varepsilon}} \frac{dJ}{J} \wedge \Phi^*\omega= \int\limits_{S^1 \times M^k } \frac{dJ}{J}\wedge\Phi^*\omega - \sum_j \int\limits_{\partial N_{j,\varepsilon}}\frac{dJ}{J}\wedge\Phi^*\omega.
\end{equation}
We need to check  that  the sum of integral over $\partial N_{j,\varepsilon}$  in the right hand side of ( \ref{E:A1})
tends to $2\pi i \int_{N_{reg}}\Phi^*\omega,$ as $\varepsilon \to 0+.$

We have locally near a point$ (\zeta_0,t_0) \in N_{reg,j}:$
$$J(\zeta,t)=(\zeta-\zeta_j(t))^m A(\zeta,t),$$
where
$$A(\zeta,t)=A_0(t)+A_1(t)(\zeta-\zeta_j(t))+\cdots$$
and $A_0(t)$  does not vanish identically. Here $t$ is taken in a neighborhood $U \subset M^k$ of $t_0.$

Then one can write
\begin{equation}\label{E:formula1}
\sum_j \int\limits_{\{|\zeta-\zeta_j(t)|= \varepsilon\} \times U}\frac{dJ}{J}\wedge
\Phi^*\omega =m \sum_j \int\limits_{\{|\zeta-\zeta_j(t)|=\varepsilon\} \times U}
\frac{d(\zeta-\zeta_j(t))}{\zeta-\zeta_j(t)}\wedge \Phi^*\omega \\
 +\sum_j\int\limits_{\{|\zeta-\zeta_j(t)|=\varepsilon\} \times U}\frac{dA}{A}\wedge
\Phi^*\omega.
\end{equation}

We are going to prove that the first term in the right hand side tends, as $\varepsilon \to 0+,$ to the right hand side of (\ref{E:8}), while the second term goes to 0. To this end,
represent the differential form $\Phi^*\omega$ as
\begin{equation}\label{E:formula}
\Phi^*\omega=\gamma dt + \sum\limits_{s=1}^k (\alpha_s d\zeta\wedge + \beta_s
d\overline\zeta) \wedge dt[s],
\end{equation}
where $dt=dt_1 \wedge \cdots dt_k$ and $[s]$ means that $dt_s$ is skipped.

Introduce the polar coordinates
$$\zeta=\zeta_j(t)+\varepsilon e^{i\varphi}.$$
Then
\begin{equation}\label{E:dzeta}
\frac{d(\zeta-\zeta_j(t)}{\zeta-\zeta_j(t)}= i \ d\varphi, \
d\zeta=d\zeta_j(t)+i\varepsilon e^{i\varphi}d\varphi, \ d\overline
\zeta=d\overline\zeta_j(t)-i\varepsilon e^{-i\varphi}d\varphi.
\end{equation}
From (\ref{E:dzeta}), $d\varphi\wedge d\zeta=d\varphi \wedge d\zeta_j(t), \
d\varphi\wedge d\overline\zeta=d\varphi\wedge d\overline\zeta_j(t)$ and hence

$$\frac{d(\zeta-\zeta_j(t))}{\zeta-\zeta_j(t)}
\wedge \Phi^*\omega= i\gamma d\varphi \wedge dt +
i d\varphi \wedge \sum\limits _{s=1}^k(\alpha_s d\zeta_j(t) +\beta_s d\overline\zeta_j(t)) dt[s],$$
where we have denoted
$$\alpha_s:=\alpha_s(\zeta_j(t)+\varepsilon e^{i\varphi},t)  \
\beta_s:=\beta_s(\zeta_j(t)+\varepsilon e^{i\varphi},t).$$

Therefore, if $\varepsilon \to 0+$ then the first term in the right hand side in (\ref{E:formula1}) tends to

\begin{equation}\label{E:Phi*}
2\pi i \ m\sum_j\int\limits_{M} \gamma(t) dt + \sum_{s=1}^k (\alpha_s(t) d\zeta_j(t) +
\beta_s(t) d\overline \zeta_j(t))\wedge  dt[s]=2\pi i \int\limits_{N_{reg,j}}\Phi^*\omega,
\end{equation}
where $\alpha_s(t)=\alpha(\zeta_j(t),t), \ \beta_s(t)=\beta(\zeta_j(t),t), \ \gamma(t)=\gamma(\zeta_j(t),t).$
Here integration first performed in a neighborhood $ t \in U$ of $(\zeta,t_0)$ belonging to a singular chain in $N_{reg,j},$ and then
the equality  extends to the integral over the entire $N_{reg,j}.$ The factor $m$ corresponds to the
 multiplicity of the chain of integration.

It remains to prove that the second term  in the right hand side of formula (\ref{E:formula1}), referring to $\frac{dA}{A},$  tends to 0  as $t \to 0+.$
It suffices to prove that for each term in the sum we have
\begin{equation}\label{to_0}
\int\limits_{|\zeta-\zeta_j(t)|=\varepsilon \times U} \frac{dA}{A}\wedge \Phi^*\omega  \to 0,
\end{equation}
as $\varepsilon \to 0+.$

The zeroth Taylor coefficient $A_0(t)$ of $A(\zeta,t)$  is real-analytic in $t$ and hence the (analytic) zero set $\{t \in U:A_0(t)=0\}$ is at most $k-1$ dimensional.
The limit, as $\varepsilon \to 0+$ is defined by the first coefficient $A_0(t)$
of the Taylor decomposition of $A(t):$
$$\lim_{\varepsilon \to 0+} \int\limits_{\{|\zeta-\zeta_j(t)|=\varepsilon\} \times U}
\frac{dA}{A}\wedge \Phi^*\omega= \lim_{\varepsilon \to 0+}
\int\limits_{\{|\zeta-\zeta_j(t)|=\varepsilon\} \times U} \frac{dA_0}{A_0}\wedge \Phi^*\omega.$$

Now we turn to the representation (\ref{E:formula}) of $\Phi^*\omega.$

First of all, notice that the integral $$\int_U \frac{dA_0}{A_0} \wedge \gamma dt=0$$
because $A$ depends only on $t$ and hence $\frac{dA_0}{A_0} \wedge dt=0$.

Now, the integrals
$$\int_U \frac{dA_0}{A_0} \alpha_s( \zeta_j(t),t))\wedge dt[s], \   \int_U \frac{dA_0}{A_0} \beta_s( \zeta_j(t),t))\wedge dt[s] $$ converge
for the following reasons.

The function $A_0(t)$ is real-analytic and hence the (analytic) zero set $\{t \in U:A(t)=0\}$ consists of pure-dimensional parts of dimensions at most $k-1.$
When the dimension is $k-1$ then the integrals converge in the sense of principal value. If the dimension does not exceed $k-2$ then the improper integrals converge as well, since the dimension of the manifold of integration $U \subset M^k$ is $k.$

We have
$$\frac{dA_0(t)}{A_0(t)} \wedge \alpha_s(\zeta,t)d\zeta_j(t) \wedge dt[s]=0,$$
 because it is a $k+1$ dimensional form in the $k-$dimensional space of parameters $t=(t_1,\cdots,t_k).$

 Therefore, (\ref{E:formula}) and the expressions for $d\zeta$ and $d\overline\zeta$ in (\ref{E:dzeta}) yield:
\begin{equation}\label{E:00}
\lim\limits_{\varepsilon \to 0+} \int\limits_{\{|\zeta-\zeta_j(t)|=\varepsilon\}\times
U}\frac{dA_0}{A_0}
\wedge \Phi^*\omega \\
=\lim\limits_{\varepsilon \to 0+} i\varepsilon\sum_s
\int\limits_{\{|\zeta-\zeta_j(t)|=\varepsilon\}\times U} \frac{dA_0}{A_0}\wedge
(\alpha_s(\zeta,t)e^{i\varphi}d\varphi -\beta_s(\zeta,t) e^{-i\varphi}d\varphi)=0,
\end{equation}
due to the convergence of the integrals of $\frac{dA_0}{A_0} \wedge \alpha_s$ and $\frac{dA_0}{A_0} \wedge \beta_s.$

Thus, finally we obtain from (\ref{E:Phi*}) and (\ref{E:00})
$$\lim_{\epsilon \to 0+}\int\limits_{\Sigma_{\varepsilon}}\frac{dJ}{J} \wedge \Phi^*\omega =2\pi  i \ \int\limits_{N_{reg}} \Phi^*\omega=2 \pi \ i \ \int\limits_{J^{-1}(0)}\Phi^*\omega.$$
Similar computations for the conjugated form $d\overline J/\overline J$ leads to the same
integral with the opposite sign which results, due to (\ref{E:J/J}), in the factor $4\pi
i $ in formula (\ref{E:8}). This finishes the proof of formula (\ref{E:8}).

\bigskip

\noindent
{\bf II. APPLICATIONS}
\section{A geometric version of Theorem \ref{T:Main}: detecting CR structure}\label{S:Lower_bounds}

In this section we present a geometric version of Theorem \ref{T:Main}. As above, $M^k$ denotes a compact closed oriented connected real analytic $k-$ manifold.

We  use the parallel terminology for families of attached analytic discs:
\begin{Definition}\label{D:degenerate_family} We say that the family $\mathcal F_{M^k}$ of the analytic discs attached to a manifold $\Lambda$ is {\bf degenerate and homologically nontrivial} if can be parametrized by a mapping
$\Phi:\Delta \times M^k \to \Lambda$ satisfying Conditions 1 and 2 of Theorem \ref{T:Main}.

\end{Definition}


Condition 1 of degeneracy in Theorem \ref{T:Main}) holds, for example, if $H_{k+1}(\Lambda)=0.$ Condition 2 of homological nontriviality geometrically means that the "`centers"'
$\Phi(0,t), \ t \in M^k$ of the attached discs $D_t=\Phi(\Delta \times \{t\})$ fill a chain  which is not the boundary of a $k+1$-chain contained in the union of the closed analytic discs $\overline D_t.$

In this section we apply Theorem \ref{T:Main} to estimating from below the CR dimension of manifolds in terms of families of
attached analytic discs.

Remind, that a manifold $\Lambda$ {\it admits } a family $\{D_t\}$ of attached analytic discs if $\Lambda=\cup_t \partial D_t.$

Applying Theorem \ref{T:Main} to a proper parametrization $\Phi$, we obtain the following geometric version, addressed to families of attached analytic discs:
\begin{Theorem}\label{T:Thm2} Let $\Lambda=\Lambda^d \subset \mathbb C^n$ be a real-analytic manifold of dimension $d \leq n.$ Suppose $\Lambda$ admits a real-analytic
$k$-dimensional  regular degenerate and homologically nontrivial family $\mathcal F=\{D_t\}$ of
attached analytic discs, parametrized by the manifold $M=M^k$ with $H_{k-2}(M)=0.$ Then the manifold $\Lambda$ has a positive CR-dimension at any point $b \in \Lambda.$
\end{Theorem}
\pf
Theorem \ref{T:Main} claims that
$rank_{\mathbb C} \ d\Phi < d$ on $\overline \Delta \times M^k.$  Pick any point $b \in \Lambda, \ b=\Phi(u), u \in \partial \Delta \times M^k$. Since $rank_{\mathbb C}\ d\Phi(u)=
\dim_{\mathbb C}(T_b\Lambda+i T_b\Lambda)$ we have
$$2d > \dim_{\mathbb R}(T_b\Lambda+iT_b\Lambda)=\dim_{\mathbb R}T_b\Lambda+\dim_{\mathbb R}(iT_b\Lambda)-dim_{\mathbb R}(T_b\Lambda \cap iT_b\Lambda)=2d-2c(b)$$
and hence  $c(b)>0.$

\qed

\begin{Remark} \label{R:after_Thm2}
One can see from the formulation of Theorem \ref{T:Main} that Theorem \ref{T:Thm2}, which, in turn,  is a corollary of Theorem \ref{T:Main},
remains true if we assume that $M$ and $\Lambda$ are compact
real-analytic chains \footnote{In the context of this article, by $k-$ chain is understood the union of $k-$ dimensional real connected oriented manifolds with assigned coefficients (multiplicities). The coefficients
appear as factors in integration over the chain of differential $k-$ forms.}
 rather than manifolds. In this case the conclusion of Theorem \ref{T:Thm2}
addresses to CR-dimension at smooth points of $\Lambda.$
\end{Remark}

In the simplest case $n=d=2, \  k=1,$ Theorem \ref{T:Thm2} implies  a characterization of compact  complex curves in $\mathbb C^2.$ Such curves, of course, must have a nonempty boundary, i.e.,
we are talking rather about subdomains of complex curves.

We will formulate this corollary in the following form:
\begin{Corollary}\label{C:Morera_curves}(Morera theorem for holomorphic curves in $\mathbb C^2$, \cite {A1},\cite {A2}).
Let $\Lambda $ be a real-analytic compact 2-manifold in $\mathbb C^2$ with $\partial\Lambda \neq \emptyset.$ Suppose that $\Lambda$ can be
covered by the boundaries $\gamma_t=\partial D_t$ of analytic discs $D_t \subset \mathbb C^2$, constituting a real-analytic regular
family depending on the parameter $t \in S^1.$
Suppose that the family $D_t$ has homologically nontrivial orbit.
Then $\Lambda$ is a one-dimensional complex
submanifold of $\mathbb C^2.$
\end{Corollary}
\pf Let us check the conditions of Theorem \ref{T:Thm2}. The condition of having homologically nontrivial orbit is just among the assumptions.

The condition of the degeneracy of the family of the analytic discs,
holds as well.
Indeed,  since $k=1$ and $d=2$ we are in the situation  $k=d-1.$  For any parametrization $\Phi$ the topological degree of the mapping $\Phi:\partial \Delta \times M^1 \to \Lambda$ is zero
because $\Lambda$ is a manifold with nonempty boundary and hence $H_2(\Lambda)=0.$

By Theorem \ref{T:Thm2}, for any $b \in \Lambda$ the tangent space $T_b(\Lambda)$ contains a complex line. Since $\dim T_b(\Lambda)=2$ the tangent space
coincides with that complex line and hence $\Lambda$ is a complex manifold. Notice that in this case $\Lambda$ contains all the attached analytic discs $D_t$ due to the uniqueness theorem for holomorphic mappings.
\qed

\begin{Remark} Since the curves $\gamma_t$ in Corollary \ref{C:Morera_curves} bound the complex manifolds, $D_t$,
 for any holomorphic 1-form $\omega$ in $\mathbb C^2$ the moment condition is fullfiled
$$\int\limits_{\gamma_t}\omega=0,$$
where $\omega$ is an arbitrary holomorphic 1-form,
and hence Corollary \ref{C:Morera_curves}  can be viewed as a Morera type characterization of complex submanifolds in $\mathbb C^2.$
\end{Remark}

\section{Examples}\label{S:Examples}

In this section we give examples showing that
all the  conditions
in Theorem \ref{T:Thm2} are essential.

In the first two examples $d=2$ and $k=1.$


\noindent
{\bf Example 1}
Let $\Lambda$ be the graph of the function $f(z)=|z|$ over the annulus $\{1 \le |z| \le 3\}:$
$$\Lambda=\{(z,|z|): 1 \le |z| \le 3\}.$$
For each $t \in \mathbb C, |t|=1,$ define
$$D_t=\{(z, 2 + Re \ t): |z| < 2+Re \ t \}.$$
The analytic discs $D_t$ are attached to the manifold $\Lambda$ and the family $\{D_t\}$ is parametrized by the parameter $t \in S^1.$ This family is degenerate since any parametrization
$\Phi:S^1 \times S^1 \to \Lambda$ has topological degree 0 due to $H_2(\Lambda)=0,$ so Condition 1 of Theorem \ref{T:Main} is fullfiled. However, Condition 2 does not hold because
any 1-cycle is contractible in the union of the discs $\cup_{t \in S^1}\overline D_t=\{(z,w)\in \mathbb C^2: Im \ w=0, 1 \le w \le 3, |z| \le w \}$ which is a truncated 3-dimensional solid cone.
So the conditions of Theorem \ref{T:Thm2} are not satisfied. The conclusion fails as well, as the truncated cone $\Lambda$ is totally real and is not a complex manifold. This example is illustrated on Figure 1, one the left.

\noindent
{\bf Example 2}
Let $\Lambda$ be the graph
$$\Lambda=\{(z,z) \in \mathbb C^2 : 1 \le |z| \le 3 \},$$
of the function $f(z)=z$ over the annulus $1 \le |z| \le 3.$
Define
$$D_t=\{(z,z) \in \mathbb C^2: |z-2t|=1 \},$$
where $t$ is the complex parameter, $|t|=1.$
The analytic discs $D_t$  are the graphs of the discs inscribed in the annulus.
As in Example 1, the 1-parameter family of analytic discs $D_t$ is parametrized by points $t \in S^1$, and each analytic disc
$D_t$ is attached (in fact, belongs) to $\Lambda$.

In this example the two-dimensional real manifold $\Lambda$ is complex, i.e. $\dim_{CR}\Lambda=1$ so that the conclusion of Theorem \ref{T:Thm2} is true.
The family $\{D_t\}$ is degenerate
because $H_2(\Lambda)=0.$ On the other hand,  the parametrization $\Phi(\zeta,t)=(\zeta+2t,\zeta+2t), |\zeta|<1, \ |t|=1$ has for $\zeta_0=0$ the  orbit $\{(2t,2t), |t|=1\}$ which bounds no 2-cycle in the union of the discs $\overline D_t.$ Therefore, the family satisfies all the conditions of Theorem \ref{T:Thm2}.
The example is illustrated on Figure 1, on the right.

\begin{figure}[h]
\centering
\scalebox{0.5}
{\includegraphics{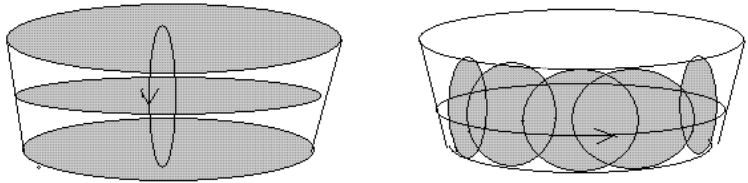}}
\label{Fig.1}
\caption{{\it  Families of attached analytic discs with homologically trivial (on the left) and homologically nontrivial (on the right) orbits}}
\end{figure}

In the next two examples $d=k=3.$

\noindent
{\bf Example 3.}
Consider the 3-manifold in $\mathbb C^3$:
$$\Lambda=\{(z_1,z_2,|z_1|^2): |z_1|^2+|z_2|^2=1\}$$
which is the graph of the function $f(z_1,z_2)=|z_1|^2$ over the unit sphere in $\mathbb C^2.$
This manifold admits 3-parameter family of attached analytic discs
$$D_t=\{(\lambda_1\zeta, \lambda_2 \zeta, |\lambda_1|^2) \in \mathbb C^3: |\zeta| <1\},$$
parametrized by the points
$$(\lambda_1,\lambda_2) \in S^3=\{(\lambda_1,\lambda_2) \in \mathbb C^2:|\lambda_1|^2 + |\lambda_2|^2=1\}$$
from the unit complex sphere in $\mathbb C^2.$

The conclusion of Theorem \ref{T:Thm2} fails, because $\Lambda$ is totally real, $\dim_{CR}\Lambda=0,$ at each point except the circle $(0,z_2,0), \ |z_2|=1.$
The condition of homological nontriviality of orbits fails for the family $\{D_t\}$ fails as well, because $\cup_{t \in S^3} \overline D_t$ is contractible and hence has the homology group $H_3=0.$

\noindent
{\bf Example 4.}
Define 3-manifold in $\mathbb C^3$ by
$$\Lambda=\{(z_1,z_2,z_1) \in \mathbb C^3: |z_1|^2+|z_2|^2=1\}.$$
The manifold $\Lambda$ is the graph of the function $f(z_1,z_2)=z_1$ over the unit sphere $S^3$ in $\mathbb C^2.$

For each point $t=(t_1,t_2) \in \frac{1}{2}S^3 \subset \mathbb C^2$ consider
the (unique) complex line $L_t$ tangent at the point $t=(t_1,t_2)$ to the 3-sphere $\frac{1}{2}S^3=\{t \in \mathbb C^2:\{|t_1|^2+|t_2|^2=\frac{1}{4}\}.$
For each point $t \in \mathbb C^2$ such that $|t_1|^2+|t_2|^2=\frac{1}{4}$ define
$$D_t:= \{z=(t_1+w_1,t_2+w_2,t_1+w_1):  |w_1|^2+|w_2|^2 \leq \frac{3}{4}\}.$$
The disc $D_t$ is the graph of the function $f(z_1,z_2)=z_1$
over the section $L_1 \cap B^2$ of the unit ball in $\mathbb C^2$ by the complex line $L_t.$ Then $\{D_t\}_{t \in S^3}$ is a 3-parametric family of analytic discs attached to $\Lambda$ and parametrized by points $t \in \frac{1}{2}S^3.$

In this example the manifold $\Lambda$ is maximally complex, i.e., $\dim_{CR}\Lambda=1,$ so the conclusion of Theorem \ref{T:Thm2} holds. All the conditions of Theorem \ref{T:Thm2} are also satisfied.
The condition of degeneracy holds for any parametrization. Indeed, $k=3,$  while $\dim \Lambda=3,$  therefore $H_{k+1}(\Lambda)=0.$
As for the second condition is concerned, choosing parametrization $\Phi$ so that $\Phi(0,t)=t$ we see that the 3-cycle $\Phi(\{0\}\times M=S^3$ is homologically nontrivial
in the spherical layer $\Phi(\overline\Delta \times M)=\{\frac{1}{4} \leq |z_1|^2+|z_2|^2 \leq 1\}$ covered by the closed discs $\overline D_t.$

\section{Lower bounds for higher CR dimensions}\label{S:Lower_general}

In the previous section we have considered the case $d \leq n$ and gave conditions for the CR-dimension of the manifold $\Lambda$ is strictly positive.
Notice that the lower bound for CR dimension is $d-n$, therefore if $d > n$ then the CR-dimension is always positive. The lower bound $d-n$ is attained in the case of generic manifolds.
Therefore the lowest nontrivial bound is $d-n+1.$

In this section we consider the case $ d \leq n$ when CR dimension, in principle, can be zero.  First, we will derive  the lowest nontrivial estimate
$\dim_{CR}\Lambda \ge d-n+1$
from  the properties of  families of attached analytic discs which $\Lambda$ admits. Then we will pass from  $d-n+1$ to higher lower bounds.

Let again $$\Phi:\Delta \times M^k \to \Lambda^d, \ \ k \geq d-1,$$ be a
real-analytic regular parametrization of the family of analytic discs $\{D_t\}_{t \in M^k}$ attached to the
real-analytic compact manifold $\Lambda^d \subset \mathbb C ^n.$

\begin{Definition} \label{D:singular} Let $\nu \leq k$ be an integer and let $C \subset M^k$ be a chain
in $M^k$ of dimension  $\dim \ C=\nu$ or  $\dim \ C =\nu-1.$
We say that that the subfamily
$$\mathcal F_{C}=\{D_t\}_{t \in C}$$
is a {\bf singular}
$\nu$-chain of attached analytic discs if it satisfies the conditions of Theorem \ref{T:Thm2}, with $M^k$ replaced by $C$,
i.e., if it is degenerate and homologically nontrivial  (see Definition \ref{D:degenerate_family}), and
$H_{\dim C -2}(C)=0.$

\end{Definition}

Remind that, as it was mentioned earlier, the meaning of the degeneracy condition in \ref{D:singular} is that either the parametrization $\Phi_{\mathcal C}$ of the family $\mathcal F_C$ decreases the dimension of $S^1 \times C$ from $\nu+1$ to $\nu$, or, in the case $\dim \ C=\nu-1$, the topological degree $\deg \Phi=0$ on $S^1 \times C.$  The second condition, of homological nontriviality, requires that the orbits $\Phi_{\mathcal C}(\zeta_0,t), \ t \in C$ represent nonzero classes in the corresponding homology group in the image $\Phi(\overline \Delta \times C).$

\begin{Definition} \label{D:tangent} Let $\mathcal F_{C}$ be a regular family  of analytic disc parametrized by a mapping $\Phi:\Delta \times C \to \mathbb C^n,$ where $C$ is a real manifold (chain). Denote
$\Lambda_C=\Phi(S^1 \times C)$.
The {\bf tangent plane} $Tan \mathcal F_{C}$ of the family $\mathcal F_C$
is the tangent plane of the manifold $\Lambda_C$  swept by the boundaries of the analytic discs:
$$Tan_b \mathcal F_{C}:=T_b\Lambda_C.$$  Here $b$ is assumed a smooth point of $\Lambda_C.$ If $C$ is a chain and $b$ is a self-intersection point of $\Lambda_C$ then we take union of tangent planes through $b.$
\end{Definition}

\begin{Definition} \label{D:passes} Let $\Pi \subset T_b\Lambda^d$ be a  linear subspace
of dimension $\nu \leq d.$
We say that the $\nu$-chain $\mathcal F_{C}$
{\bf passes through the point $b$ in the  direction $\Pi$} if $\Pi \subset Tan_b \mathcal F_{C}.$
\end{Definition}


\begin{Definition} \label{D:admissible} We will call a real $\nu$-plane $\Pi \subset T_b\Lambda$  {\bf admissible} if $2n-d+\nu$ is even and there exists a
complex linear subspace $P \subset T_b \mathbb C^n,$ of real dimension $\dim_{\mathbb R}P =2n-d+\nu,$ such that $\Pi=P \cap T_b\Lambda.$
\end{Definition}

The following theorem applies to the case $d \ge n$ and gives conditions to ensure that the CR-dimension is  bigger than the generic
dimension $d-n.$

\begin{Theorem}\label{T:Thm3} Let $ d \ge  n$ and let $\Lambda=\Lambda^d \subset \mathbb C^n$ be a real-analytic  manifold. Suppose that $\Lambda$
admits a  real-analytic regular $k$-parametric, $k \ge d$, family $\mathcal F=\mathcal F_{M^k}$ of attached analytic discs.
Suppose that the family $\mathcal F$ has the following property:
for any point  $b \in \Lambda$  and for any
$(2n-d)$- dimensional admissible plane $\Pi \subset T_b(\Lambda)$
there exists a singular $(2n-d)$-dimensional chain
 $\mathcal F_{C} \subset \mathcal F$  passing through $b$
in the direction $\Pi$. Then at any point $b \in \Lambda$ the CR-dimension satisfies
$$c(b)=\dim_{\mathbb C} T^{\mathbb C}_b \Lambda \ge d-n+1.$$
In other words, $\Lambda$ is nowhere generic.
\end{Theorem}
The proof of Theorem \ref{T:Thm3} follows from Theorem \ref{T:Thm2}:

\pf

In this case the parameters in Definition \ref{D:admissible} are: $\nu=2n-d$ and $\dim_{\mathbb R}=2n-d+\nu=2(2n-d).$

Suppose that $\Lambda$ is generic at  a point $b \in \Lambda,$ i.e.,  $c(b)=d-n.$

The $d$-dimensional tangent plane
decomposes into the direct sum of a complex plane of the complex dimension $d-n$ and a totally real (free of nonzero complex subspaces)
plane of the real dimension $d-2(d-n)=2n-d$:
\begin{equation}\label{E:Tb}
T_b \Lambda= \Pi_{\mathbb C}^{d-n} \oplus \Pi_{\mathbb R}^{2n-d}.
\end{equation}
Let us show that the second summand is an admissible linear subspace.

First of all, since $\Pi_{\mathbb C}^{d-n}$ is the maximal complex linear subspace contained in $T_b\Lambda$, we have
$$i\Pi_{\mathbb R}^{2n-d} \cap T_b\Lambda=0.$$ Indeed, if $p \in \Pi_{\mathbb R}^{2n-d}$ and $ip \in T_b\Lambda$
then  the complex line $\mathbb C \cdot p$ belongs to  $T_b\Lambda$ and $\Pi_{\mathbb C}^{d-n} \oplus \mathbb C \cdot p$ is a complex subspace in $T_b\Lambda$ contained $\Pi_{\mathbb C}^{d-n}$. Therefore, due to maximality,
$ p \in \Pi_{\mathbb C}^{d-n}$ and then $p=0$ since the decomposition (\ref{E:Tb}) is direct.

Now consider the complexified linear space
$$P:=\Pi_{\mathbb R}^{2n-d} \oplus i\Pi_{\mathbb R}^{2n-d}.$$
Since the $\nu$-dimensional,  $\nu=2n-d$,  linear space $\Pi_{\mathbb R}^{2n-d}$ is totally real, we have
$$\Pi_{\mathbb R}^{2n-d} \cap i\Pi_{\mathbb R}^{2n-d}=0.$$
Then
$$\dim P=dim_{\mathbb R}\Pi^{2n-d}+dim_{\mathbb R}(i\Pi^{2n-d}=2(2n-d)=2n-d+\nu,$$ so the dimension of $P$ is as required in Definition \ref{D:admissible}.
If $p \in P \cap T_b\Lambda$ then $p=q+ir,$ where $ \ q,r \in \Pi_{\mathbb R}^{2n-d}$ and since $p,q \in T_b\Lambda$ we have
$ir=p-q \in T_b\Lambda.$ But $ir \in i\Pi_{\mathbb R}^{2n-d}$ and hence $ir \in i\Pi_{\mathbb R}^{2n-d} \cap T_b\Lambda=0.$ Thus, $r=0$ and hence $p=q \in \Pi_{\mathbb R}^{2n-d}.$ This shows that
$$P \cap T_b\Lambda \subset \Pi_{\mathbb R}^{2n-d}.$$
The converse inclusion is obvious. Therefore $$\Pi_{\mathbb R}^{2n-d}=P \cap T_b\Lambda$$ is an admissible linear space.

Now, using the condition, we conclude that  there exists a singular $\nu$-chain $\mathcal F_{C}, \ \nu=2n-d,$  passing through the point $b$ in  the direction
$\Pi^{2n-d}.$
Consider the manifold
$$\Lambda^{2n-d}:=\Phi(S^1 \times C)$$
and apply Theorem \ref{T:Thm2} to the family
$\mathcal F_{C}$ of the analytic discs attached to $\Lambda^{2n-d}.$

First of all, observe that for the manifold $\Lambda^{2n-d}$ we are in the range
of dimensions required in Theorem \ref{T:Thm2}, because $d \geq n$ and hence the real dimension of $\Lambda^{2n-d}$ is less than
the complex dimension of the ambient space: $2n-d \leq n$. All other conditions of Theorem \ref{T:Thm2} are satisfied for the manifold $\Lambda^{2n-d}$ and the family $\mathcal F_{C}.$

By Theorem \ref{T:Thm2} the CR-dimension of $\Lambda^{2n-d}$ satisfies the inequality $c_{\Lambda^{2n-d}}(b) >0.$
In other words, the tangent plane $T_b\Lambda^{2n-d}$ contains
a complex line. But this tangent plane is contained in the totally real $(2n-d)$-plane
$ \Pi_{\mathbb R}^{2n-d}$ which is
free of nontrivial complex subspaces. This contradiction with the assumption $c(b)=d-n$ shows that $c(b)\ge d-n+1.$

\qed

\begin{Remark}

\begin{enumerate}
\item In the case $d=n$ Theorems \ref{T:Thm2} and \ref{T:Thm3} coincide. Indeed, then
$2n-d=d$  and therefore the only $(2n-d)$- direction is the whole tangent plane $T_b \Lambda.$
\item Due to real analyticity, it suffices to require that the condition in Theorem \ref{T:Thm3} holds for the points $b$  from an open set $U \subset \Lambda.$
\end{enumerate}
\end{Remark}

The next result generalizes Theorem \ref{T:Thm3} for higher CR-dimensions.
It gives conditions for the CR-dimension to be at least $q,$ where $q \geq d-n$ is a given natural number, $q \le n.$

\begin{Theorem}\label{T:Thm4}
Let $d \geq n.$ Let the manifolds $\Lambda=\Lambda^d$ and $M^k$, and the family
$\mathcal F_{M^k}$  be as in Theorem \ref{T:Thm3}. Fix natural
$q, \ d-n \leq q \le n.$  Suppose that for any $b \in \Lambda$
and for any admissible real plane $\Pi \subset T_b\Lambda$ of the dimension $\nu:=d-2q+2,$
there exists a singular $\nu$-chain  $\mathcal F_{C} \subset \mathcal F$
passing through $b$ in the direction $\Pi.$ Then  for any $b \in \Lambda$ holds:
$$ c(b) \ge q, \  b \in \Lambda.$$
\end{Theorem}

\pf
Let $b \in \Lambda$. Let $c(b)=\dim_{CR} \Lambda(b)$ be the CR-dimension at the point $b$. We know that always $c(b) \geq d-n$.
As in the proof of Theorem \ref{T:Thm3}, the tangent plane decomposes in the direct sum
$$ T_b\Lambda=\Pi^{c(b)}_{\mathbb C} \oplus \Pi_{\mathbb R}^{d-2c(b)}$$ of
a complex subspace and a totally real admissible subspace.

Now suppose that, contrary to the assertion, $$c(b) \leq q-1.$$
Then we have $$d-2c(b) \geq d- 2(q-1)=\nu.$$

Take an arbitrary admissible $\nu-$ plane $\Pi^{\nu} \subset \Pi_{\mathbb R}^{d-2c(b)}.$
Applying, as in the previous proof, Theorem \ref{T:Thm2} to the $\nu-$ chain $\mathcal F_C,$ passing through
$b$ in the direction  $\Pi^{\nu},$ we conclude that  $\Pi^{\nu}$ must contain a nonzero complex subspace.
However, this is impossible because $\Pi_{\mathbb R}^{d-2c(b)}$ is totally real. This contradiction shows that $c(b) \ge q.$

\qed

\section{Characterization of holomorphic manifolds and their boundaries}\label{Ss:complex_boundaries}

Now we  turn  to characterization of  complex manifolds or boundaries of complex submanifolds of the ambient space $\mathbb C^n.$

If $\Lambda^d$ is an orientable  submanifold of $\mathbb C^n$ of even real dimension $d=2p$
then the CR dimension $q$ of $\Lambda$ is at most $p$ and  $q=p$ if and only if $\Lambda^d$ is holomorphic (complex) submanifold in $\mathbb C^n.$

When $d$  is odd, $d=2p-1,$ and $d>1$ then
the CR dimension $q$ is at most $p-1$ , and due to theorem of Harvey and Lawson \cite{HL} the CR dimension is maximal possible, $q=\dim_{CR}\Lambda^d=p-1,$ if and only if $\Lambda=\partial V^p$ where $V^p$ is a complex manifold in $\mathbb C^n,$ of complex dimension $p.$

In the case $d=1$, i.e., $\Lambda^d$ is a curve, according to the result of Wermer \cite{W1,W2},  $\Lambda=\partial V^1$ if and only if the moment condition holds: $\int_{\Lambda}\omega=0$ for all holomorphic 1-forms $\omega.$
The above results hold not only for manifolds but for chains as well.

Thus, when $\dim \Lambda=d >1$ then for  both even and odd dimensions $d$  holomorphic manifolds or their boundaries correspond to maximally complex manifolds, i.e. manifolds with maximally possible CR dimension.

Let us start with the characterization of holomorphic manifolds ($d$ even).
Take $q=d=2p.$  Then applying  Theorem \ref{T:Thm4} with  $\nu=d-2q+2=2,$ we arrive at
\begin{Theorem}\label{T:Thm5} Let $\Lambda=\Lambda^{2p}$ be a real-analytic orientable manifold, admitting a real-analytic regular family $\mathcal F$ of attached analytic discs. Suppose that for any $b \in \Lambda$ and for any admissible tangent 2-plane $\Pi$ there exists a singular 2-chain $\mathcal F_C \subset \mathcal F$, passing through $b$ at the direction $\Pi$. Then $\Lambda$ is a holomorphic manifold of the complex dimension $p.$
\end{Theorem}
Remind, that the admissible 2-planes are intersections of the tangent spaces $T_b(\Lambda)$ with complex $(n-p+1)$-planes.
Theorem \ref{T:Thm5} generalizes Corollary \ref{C:Morera_curves} of Theorem \ref{T:Thm2} from the dimension $n=2$ of the ambient space to arbitrary dimension $n.$

Now we turn  to the characterization of the boundaries of holomorphic manifolds or chains. This case corresponds to odd dimensions
$\dim \Lambda=d=2p-1.$ We set $q=p-1.$ Then the parameter $\nu$ in Theorem \ref{T:Thm4} becomes
$\nu=d-2q+2=(2p-1)-2(p-1)+2=3.$
The admissible 3-planes are obtained as intersections  of $\Lambda$ with complex $(n-p+2)$-planes.

According to the cited above result of Harvey and Lawson \cite{HL}, if $dim_{CR}\Lambda^{2p-1}=p-1,$  then there exists $p$-dimensional complex chain $V \subset \mathbb C^n$ such that
$\Lambda^{2p-1}=\partial V.$ Applying Theorem \ref{T:Thm4} we obtain

\begin{Theorem} \label{T:Thm6} Let $\Lambda^{2p-1}, p>1,$ be a closed orientable real-analytic manifold in $\mathbb C^n$ admitting a real-analytic regular family $\{D_t\}_{t \in M^k}$
of  attached analytic discs. Suppose that for any admissible tangent 3-plane $\Pi \subset T_b(\Lambda)$ there exists a singular
3-chain $\mathcal F_{C}$ passing  through $b$ in the direction $\Pi.$
Then $\Lambda$ bounds a holomorphic $p-$ chain $V^p$ , $\Lambda^{2p-1}=\partial V^p.$
\end{Theorem}

Theorem \ref{T:Thm6} takes most simple form in the case $p=2.$ Let us formulate the corresponding statement in a rather explicit form.
\begin{Theorem} (Theorem \ref{T:Thm6} for the case $p=2$). Let $\Lambda^3 \subset \mathbb C^n$ be a real-analytic closed orientable 3-manifold. Suppose that there is a
real analytic family of analytic discs $D_t$ which can be parametrized by a 3-dimensional real-analytic closed manifold $M^3$ by means of the real-analytic mapping $\Phi:\overline \Delta \times M^3 \to \mathbb C^n$ which is regular with respect to to $t \in M^3.$ Suppose that
\begin{enumerate}
\item the curves $\gamma_t=\partial D_t$ are contained in $\Lambda ^3 $ and cover it, i.e. $\Lambda^3=\cup\gamma_t, \ t \in M^3.$
\item  for some fixed $\zeta_0 \in \Delta$ the cycle $\Phi(\{\zeta_0\} \times M^3)$ bounds no 4-chain in $\cup \overline D_t, t \in M^3$ (i.e., the family $D_t$ has a homologically nontrivial orbit).
\end{enumerate}
Then there is a holomorphic chain $V^2$ of complex dimension 2 such that $\partial V^2=\Lambda^3.$
\end{Theorem}
The condition 2 of having nontrivial orbit is crucial. The corresponding example is essentially given in Section \ref{S:Examples} (Example 3).
In the context of Theorem \ref{T:Thm6},
Example 3 illustrates the following.

When $\Lambda ^3$ is the graph of a function over the boundary $\partial \Omega$ of a complex domain $\Omega$
then the condition of $\Lambda^3$ to be the boundary of a holomorphic chain recasts in the condition of $f$ to be the boundary value of a holomorphic function in $\Omega.$
Thus, the counterexample can be delivered by a function that extends holomorphically in a family of analytic discs, but has no global holomorphic extension in the domain.

Such a function can be taken, for example, $f(z_1,z_2)=|z_1|^2$ on the boundary of the unit ball $\Omega=B^2$ in $\mathbb C^2.$ Then the graph $\Lambda^3=graph_{\partial B^2}f \subset \mathbb C^3$
is the boundary of no holomorphic
2-chain $V^2$, because $f$ does not possess holomorphic extension from the complex sphere $\partial B^2$ to the unit ball $B^2.$
The function $f(z_1,z_2)$ is constant on any circle $C_t=\{(\zeta a, \zeta b): |\zeta|=1\}, \ t=(a,b) \in \partial B^2$ and hence $\Lambda^3$ is covered by the
boundaries of analytic discs $D_t=\{(\zeta a, \zeta b, |a|^2): |\zeta|<1\}.$ The condition of homological nontriviality is not fullfiled because
the union $\cup\overline D_t, \ t \in \partial B^2$ is contractible and hence any cycle in it is trivial.

Dolbeault and Henkin \cite{DH0}, \cite{DH} proved that the
closed maximally complex manifold $\Lambda=\Lambda^{2p-1} \subset \mathbb CP^n $
bounds a holomorphic chain,
if for any complex $(n-p+1)$-plane $P,$ transversally intersecting $\Lambda,$
the (one-dimensional) intersection $\Lambda \cap P$ ,  bounds a complex 1-chain.  T.-C.Dinh \cite{D1},\cite{D2},\cite{D3}
refined the result, by eliminating the condition for $\Lambda$ to be maximally complex and by using a
narrower family of the complex $(n-p+1)$-sections.
This result, in the case when the manifold $\Lambda$ is contained  in $\mathbb C^n$,
corresponds to  a special case of Theorem \ref{T:Thm5}, for "`linear"'  attached analytic discs which are bounded by intersections of the manifold $\Lambda$ with linear complex spaces.
The required by Theorem \ref{T:Thm5} singular chains $\mathcal F_{C},$ passing through an arbitrary point $b,$ correspond
to complex $(n-p+1)$-subspaces of the complex $(n-p+2)$-subspace defined a given admissible direction. Our condition of  homological nontriviality appears in  \cite {D1},\cite{D2}
as the condition that the analytic discs (complex sections) are disjoint from a fixed compact $(n-p+1)$-linearly convex set $Y.$

\begin{Remark}

\begin{enumerate}
\item
Due to real analyticity, it suffices to require in Theorems \ref{T:Thm3}-\ref{T:Thm6}
that the singular chains condition hold for the points $b$ from an open set $U \subset \Lambda.$
\item
Most simply Theorem \ref{T:Thm6} looks in the case $n=3, p=2, d=2p-1=3$. In this case
Theorem \ref{T:Thm6} says that a sufficient condition for the
real-analytic closed 3-dimensional surface $\Lambda$ to be a boundary $\Lambda=\partial V$ of a complex 2-chain $V$
is that $\Lambda$ admits 3-parameter regular real-analytic family of attached analytic discs with homologically nontrivial orbit. The condition of degeneracy holds because in this case $k+1=3+1 > d=3$. Notice that, vice versa,
if $\Lambda$ is a boundary of complex 2-chain $V$ then the homologically nontrivial family of attached analytic discs exists; it can be constructed from the analytic discs in $V,$ attached to the boundary  $\partial V=\Lambda$
and located nearby $\partial V.$
\end{enumerate}
\end{Remark}
\section{Morera type theorems for CR-functions}\label{S:applications_Morera}

In this section we apply the above results about CR-dimensions of manifolds to the case, when the manifold $\Lambda$ is the graph of a function. Then the results for manifolds, translated to the language of functions,
lead to  characterization of CR-functions
in terms of analytic extendibility into attached analytic discs, or, equivalently, in terms of zero complex moments on the boundaries of the attached discs.
Everywhere, the parametrizing manifold $M$ is assumed satisfying the conditions from Theorem \ref{T:Main}.

\subsection {The case of CR dimension one}

We start with the case of generic manifolds of real dimension $d=n+1,$ where $n$ is the complex dimension of the ambient space $\mathbb C^n.$ In this case, the CR-dimension equals $(n+1)-n=1.$
The following theorem follows from  Theorem \ref{T:Thm2}, applied to real manifolds which are graphs.
\begin{Theorem}\label{T:Thm7} Let $\Omega=\Omega^d \subset X^n$ be a real-analytic generic submanifold of the real dimension $d=n+1.$ Let
$\{\Omega_t\}_{t \in M^k}$ be a real-analytic family of attached analytic discs, degenerate and having homologically nontrivial orbit (Definition \ref{D:degenerate_family}),
parametrized by a compact real-analytic connected closed $k$-manifold $M^k.$
Let $f$ be a real-analytic function on $\Omega$ such that the restriction
$f\vert_{\partial \Omega_t}$ analytically extends in the analytic disc
$\Omega_t$, for any $ t \in M,$ or, equivalently,
$$\int\limits_{\partial \Omega_t} f \omega=0$$ for any holomorphic 1-form in $\mathbb C^.$
Then $f$ is CR-function on $\Omega$, i.e. $f$ satisfies the boundary CR-equation everywhere on $\Omega.$
\end{Theorem}

\pf
Let $\Lambda$ be the graph of the function $f$,
$$\Lambda=graph_{\Omega}f=\{(z,f(z): z \in \Omega\} .$$
Then $\Lambda$ is a real-analytic submanifold of the $(n+1)-$ dimensional complex space $\mathbb C ^n \times \mathbb C.$ Thus, we are in the situation of Theorem \ref{T:Thm2} because the real dimension $n+1$ of the manifold $\Lambda$ equals to the complex dimension of the ambient complex space.

Since $\Omega$ is generic, we have $\dim_{CR}\Omega=d-n=1.$
As the CR-dimension of $\Lambda$ is concerned, the guaranteed lower bound is trivial:
$$c_{\Lambda}(b) \geq \dim_{\mathbb R}\Lambda -(n+1)= (n+1)-(n+1)=0, \ b \in \Lambda.$$  However, applying Theorem \ref{T:Thm2},
we will show that in fact $c_{\Lambda}(b)$ is strictly positive at any point $b \in \Lambda.$

Denote $Q_f$ the lifting mapping
$$Q_f:\Omega \to \Lambda, \ \ Q_f(z)=(z,f(z)).$$
Let $\Psi(\zeta,t), \zeta \in \Delta, t \in M^k,$ be the parametrization of the family $\{\Omega_t\}.$

If $F_t$ is the analytic extension of the function $f$ into the analytic disc $\Omega_t$ then the composition mapping
$$\Phi_t = F_t \circ \Psi_t$$
defines the parametrization of the family of analytic discs $$D_t=(F_t \circ \Psi(\cdot,t))(\Delta)$$
attached to the manifold $\Lambda.$ There analytic  discs $D_t$ are the graphs
$$D_t=graph_{\Omega_t} F_t$$
of the analytic extensions $F_t$ into analytic discs $\Omega_t.$

It can be readily checked that the conditions for the manifold $\Omega$ and for the parametrization
$\Psi$ in Theorem \ref{T:Thm7} translates as the corresponding conditions in Theorem \ref{T:Thm2}, for the
manifold $\Lambda$ and the parametrization $\Phi$ (with the dimension $n$ of the ambient space replaced by $n+1.$)

Then Theorem \ref{T:Thm2} implies the estimate
$$\dim_{CR} \Lambda \geq 1.$$
This means that at any point $u \in \Omega$  the differential $dQ_{f}(u)$ does not decrease the CR-dimension and
maps the one-dimensional complex subspace in $T_u \Omega$ to a one-dimensional complex subspace of $T_{(u,f(u))} \Lambda$. In other words,
the differential $df(u)$ is a complex linear map on $T_u^{\mathbb C}\Lambda$ and hence $f$ satisfies the tangential CR-equation at the point $u.$
\qed

\subsection{Special case of Theorem \ref{T:Thm7}: the strip-problem}\label{S:the_strip_problem}
The simplest special case $n=1$  of Theorem \ref{T:Thm7} gives an answer, for real-analytic case, to a question, known as the {\it strip-problem} and discussed in subsection \ref{S:the_strip_problem}.
The detailed proof is given in the article \cite {A1}. We will give here a brief version of the proof, referring the reader to \cite{A1} for more details and the references. Since we restrict ourselves
by  closed families of analytic discs,  in the case of one-parameter families the parametrizing manifold can be taken the circle $S^1.$

Now we have $n=1, k=1, d=2$ and the manifold $\Lambda$ is a domain in the complex plane.
\begin{Corollary} \label{C:4.2} Let $\Omega$ be a compact domain in the complex plane, covered by 1-parameter regular real-analytic family
of Jordan curves $\gamma_t, \ t \in S^1,$ such that no point in $\Omega$ is surrounded by all the curves $\gamma_t$. Let $f$ be a real analytic function
on $\Omega$ such that all complex moments
$$\int\limits_{\gamma_t} f(z)z^m dz=0, \ \  \forall t , k=0,1,\cdots.$$
Then $f$ is holomorphic in $\Omega.$
\end{Corollary}

\begin{Remark}
\begin{enumerate}
\item The condition of regularity of the family $\gamma_t$ means in our case that the velocity vector is not proportional to the tangent vector in all point in $\Omega$ except for the boundary ones.
Parametrize the family by real-analytic family of conformal mappings $\Phi_t:\Delta \to D_t, \ \gamma_t=\partial D_t.$ Then the regularity means
$$Im \frac{\partial_t\Phi(e^{i\psi},t)}{\partial_{\psi}\Phi(e^{i\psi},t)} \neq 0$$
for all $\psi$ and $t$ such that $\Phi(e^{i\psi},t) \notin \partial \Omega.$
\item
Vanishing complex moments of $f$ on the curves $\gamma_t$ is equivalent to analytic extendibility of $f$ inside the domain $\Omega_t$ bounded by $\gamma_t.$
The condition for the family of curve $\gamma_t$ means  that $cap_{t \in S^1} \overline \Omega_t=\emptyset.$
\end{enumerate}
\end{Remark}
\pf {\bf of Corollary \ref{C:4.2}}
Let $\Psi$ be a regular parametrization of the family $\Omega_t.$
As above, we can take $\Psi_t=\Psi(\cdot,t)$- the conformal mapping of
the unit disc $\Delta$ onto $\Omega_t$. These mappings are assumed to be chosen real-analytically depending on the parameter $t$.

Let us check the conditions of Theorem \ref{T:Thm7}.
First of all,  since the compact domain $\Omega$
lies in $\mathbb R^2$ and has nonempty boundary, the Brouwer degree of the mapping $\Psi: S^1 \times M^1 \to \Omega \subset \mathbb R^2$ equals 0. Thus, the parametrization $\Phi$ is homologically
degenerate (Definition \ref{D:degen}).

Now, CR-functions on $\Omega$ are simply holomorphic functions.
To derive Corollary \ref{C:4.2} from Theorem \ref{T:Thm7}, it  only remains to check that
the family  $\Omega_t$ has homologically nontrivial orbit. It is rigorously done in \cite{A2} (the families with homologically nontrivial orbits are called there homologically nontrivial).
We  will give here the sketch of the argument.

Homological triviality of the family of the domains $\Omega_t$ means that the  1-cycle
$$c=\Psi (\{0\} \times M) \subset \tilde \Omega:=\cup_{t \in M^1} \overline \Omega_t$$
is homologically trivial in $\tilde \Omega.$
Then $c$ can be contracted, within the domain $\tilde \Omega$, to a point.

Denote $b$ the point to which $c$ contracts.
By lemma about covering homotopy (cf.\cite {Hu},pp.61-66,) the preimage cycle $C=\{0\} \times M$ can be correspondingly deformed to a nontrivial cycle $C^{\prime} \subset \overline\Delta \times M^1.$ The cycle $C^{\prime}$
necessarily meets each closed disc $\overline \Delta \times \{t\}$. On the other hand,
$\Psi(C^{\prime})=\{b\}.$ Therefore $b$ belongs to each domain $\overline \Omega_t=\Psi(\overline \Delta \times \{t\}).$  We have obtained  a contradiction with
the condition that the closed domains $\overline D_t$ have empty intersection. This proves that our family has homologically nontrivial orbit.
Thus, the family $\Omega_t$ satisfies all the conditions of Theorem \ref{T:Thm7} and therefore Corollary \ref{C:4.2} follows.
\qed

\subsection{The case of arbitrary CR dimension.}

The following theorem is a version of Theorem \ref{T:Thm4}  for the case when the manifold is a graph:
\begin{Theorem} \label{T:Thm8} Let $\Omega=\Omega^d \subset \mathbb C^n$ be a real-analytic CR-manifold, $\dim_{CR} \Omega=q.$
Suppose that $\Omega$ is covered by the boundaries of
a regular real-analytic family $\mathcal F=\{\Omega_t\}_{t \in M^k}, \ k \geq d-1,$ of attached analytic discs. Suppose that
for any $b \in \Omega$ and for any admissible real linear space $\Pi=\Pi^{d-2q+2} \subset T_b\Omega$ there exists a singular $(d-2q+2)$-chain
$\mathcal F_{C} \subset \mathcal F$ passing through $b$ in the direction $\Phi$
(i.e. satisfying the conditions of Theorem \ref{T:Thm4} with the parameter $q.$)
Let $f$ be a real-analytic function on $\Omega$ and suppose that $f$ satisfies the Morera condition
$$\int\limits_{\partial \Omega_t} f \omega=0, \ \forall t \in M^k,$$
for arbitrary holomorphic 1-form $\omega$,  or, equivalently,
that $f$ analytically extends in each analytic disc $\Omega_t.$ Then $f$ is CR-function on $\Omega.$
\end{Theorem}
\pf
Denote $$\Lambda=\Lambda^d=graph_{\Omega} f \subset X^n \times \mathbb C.$$ Let $F_t$ be the analytic extension of $f$ into $D_t.$ Then the graphs
$$D_t=graph_{\Omega_t} F_t \subset X^n \times \mathbb C$$
are analytic discs and $\Lambda$ is covered by their boundaries.

The manifolds
$\Omega$ and $\Lambda$ are linked via the lifting diffeomorphism $Q_f:\Omega \to \Lambda$
given by
$$ Q_f(u)=(u,f(u)), \ u \in \Omega.$$
Analogously, the analytic disc $D_t$ is the image of $\Omega_t$ under the lifting diffeomorphism
$$Q_{F_t}(u)=(u,F_t(u)).$$
Since $F_t=f$ on $\partial \Omega_t$, we have $Q_f\vert_{\partial\Omega_t}=Q_{F_t}\vert_{\partial \Omega_t}$  and therefore
$$\partial D_t=\partial Q_{F_t}(\Omega_t)=Q_{F_t}(\partial \Omega_t)=Q_f(\partial \Omega_t) \subset \Lambda, $$
which means that the analytic discs $D_t$ are attached to $\Lambda.$

One can readily check that  the conditions of Theorem \ref{T:Thm4}, for the manifold $\Omega^d$ and for the family $\{\Omega_t\}_{t \in M^k}$ ,
imply same type conditions of Theorem \ref{T:Thm4}, for  the graph $\Lambda^d$ and  for the family  $\{D\}_{t \in M^k}$,  with the same parameter $q.$
Applying Theorem \ref{T:Thm4} one concludes that the CR dimension of $\Lambda$ satisfies
$$c_{\Lambda}(b) \geq q.$$
But this means that the mapping $Q_f$ does not decrease CR dimension and hence $Q_f$ is CR-mapping. Then  $f$ is a CR-function
as the superposition $f=\pi_2 \circ Q_f$
of two CR-mappings: $Q_f$ and the projection $\pi_2(u,w)=w$.
\qed
\begin{Remark} Most simply Theorem \ref{T:Thm8} formulates for the case of real hypersurfaces $\Omega$ in $\mathbb C^2.$ In this case $d=3, q=1, d-2q+2=3$, so that the family $\{\Omega_t\}$ is  3-dimensional.
The conditions of Theorem now addresses to the entire family $\{\Omega_t\}$  and says that it should be degenerate and homologically nontrivial.
\end{Remark}
\subsection{Special cases of Corollary \ref{C:V^p}: $n-$dimensional strip-problem and Globevnik-Stout conjecture}
In the special case $d=2q=2n,$ the manifold $\Omega$ is a domain in $\mathbb C^n$ and we obtain from Theorem \ref{T:Thm8} a test of holomorphicity, which generalizes Corollary \ref{C:4.2} from $\mathbb C$  to $\mathbb C^n.$
In this special case $d-2q+2=2,$ so that the admissible 2-directions are complex lines and  we have

\begin{Corollary}($n$-dimensional strip-problem) \label{C:n-dim-strip}.
Let $\Omega$ be a domain in $\mathbb C^n$ covered by the boundaries $\partial \Omega_t$ of the analytic discs, constituting
a regular real-analytic family $\mathcal F_M.$ Suppose that $\mathcal F_M$ contains
singular 2-chains $\mathcal F_C \subset \mathcal F$ passing through each point $b \in \Lambda$ (due to real-analyticity, $b$ can be taken from an open set) in any prescribed one-dimensional complex direction.
Let $f$ be a real-analytic function in $\overline \Omega$ and assume that
$$\int_{\partial \Omega_t} f\omega=0$$
 for every $t \in M$ and every holomorphic 1-form.
Then $f$ is holomorphic in  $\Omega.$
\end{Corollary}

Another interesting
special case of Theorem \ref{T:Thm8} is when $\Omega$ is the boundary of a complex manifold. In this case  $d=2p-1, \ q=p-1,
\ d-2q+2= 3$.
By Bochner-Hartogs theorem \cite{Hor}, the smooth boundary values of holomorphic functions on the (smooth) boundary of a domain in $\mathbb C^n$ coincide with CR-functions and therefore we obtain:

\begin{Corollary}\label{C:V^p} Let $V \subset  \mathbb C^n$ be a $p$-dimensional complex manifold with the real-analytic boundary $\partial V=\Omega.$
Suppose that $\Omega$ admits a real-analytic regular family $\{D_t\}_{t \in M^k}$ of attached analytic discs and assume that
there are 3-dimensional singular chains in $\{D_t\}$ passing through $\Omega$ in any admissible 3-dimensional direction.
If $f$ is a real-analytic function on $\Omega$ and $f$ admits analytic extension into each analytic disc $D_t, \ t \in M^k$,
then $f$ extends from $\Omega$ as a holomorphic function in $V.$
\end{Corollary}

A few comments. First, notice that all the analytic discs $D_t$ attached to $\partial V$ are necessarily contained in $V.$  Homological nontriviality, which is one of the conditions for singular chains,
can be provided, for example, by demanding  that all  the  closed analytic discs
$D_t, \ t \in M^k,$  fill  $V \setminus V_0$ where $V_0 \subset V $ is an open subset ( a "hole"). In Corollary \ref{C:GSC}, which is a special case of Corollary \ref{C:V^p} and is presented below,  the family is of exactly that type.

First result for such type of families was proved by Nagel and Rudin \cite{NR}. They showed that if $f$ is a continuous function on the unit sphere
in $\mathbb C^n$ having analytic extension into any complex line on the fixed distance to the origin then $f$ is the boundary value of
a holomorphic function in the unit ball. The proof essentially used harmonic analysis in the unitary group and did not extend to
non-group invariant case.

Globevnik and Stout formulated in \cite{GS} the problem of generalization of Nagel-Rudin theorem for arbitrary domains.
T.-C.Dinh \cite {D2} proved Globevnik-Stout conjecture under assumption of non real-analyticity of the line sections.
Baracco, Tumanov and Zampieri \cite{BTZ} confirmed the conjecture  with tangent Kobayashi geodesics in the place of linear sections.
(see \cite{A2} for the extended references).

Now we will present a result for general case, even in a form stronger than Globevnik-Stout conjecture, however under assumption of real-analyticity, as everywhere is this article.
So, the special case of Corollary \label{C:V_p} when $V$ is a domain in $\mathbb C^n,$ is the following statement
(see \cite{A1}, \cite{A2}):

\begin{Corollary} (a strong version of Globevnik-Stout conjecture for real-analytic case) \label{C:GSC}
Let $D \subset \mathbb C^n$ be a bounded domain with real-analytic boundary $\partial D.$ Suppose that $\partial D$ is covered by the boundaries of attached analytic discs $\{D_t\}_{t \in M},$ constituting a real-analytic regular family containing 3-dimensional singular 3-chains passing through $\partial D$ in any prescribed admissible 3-dimensional direction. If $f \in C^{\omega}(\partial D)$ extends analytically in each disc $D_t$ then
$f$ is CR function and hence is the boundary value of a holomorphic function in $D.$
\end{Corollary}

Globevnik-Stout conjecture is a special case of Corollary \ref{C:GSC} when $D$ is convex and the discs $D_t$ are cross-sections of $D$ by complex lines tangent to the boundary of a fixed convex real-analytic
subdomain $D^{\prime} \subset D:$
\begin{Corollary}\label{C:GS_lines} Let $D^{\prime} \subset D$ be two strictly convex bounded domains with real analytic boundaries. If $f \in C^{\omega}(\partial D)$ admits analytic extension in each
section $L \cap D$ by complex line $L$ tangent to $\partial D^{\prime}$ then $f$ continuously extends in $\overline D$ as a function holomorphic in $D.$
\end{Corollary}
\pf

Consider first the case $n=2.$
The family of complex lines $L_t$ tangent to $\partial D^{\prime}$ (at the point $t \in \partial D^{\prime}$) is regular and real-analytic.
This family is 3-parametric, therefore we are in the situation $k+1=4>d=3$ and hence the family is homologically degenerate on the boundary.
Moreover, it is homologically nontrivial, because $\partial D^{\prime}$ is the (3-dimensional) orbit of the family, which
bounds no domain in the union of the analytic discs $\cup(L_t \cap D)=\overline D \setminus D^{\prime}$ because of the hole $D^{\prime}.$ By Corollary \ref{C:GSC} $f$ extends to $D$ as a holomorphic function. In particular, $f$ satisfies tangential CR equation on $\partial D.$

Let $n>2.$ We apply the just proved statement to any section $\Pi \cap D$ by complex 2-plane $\Pi \subset \mathbb C^n,$ intersecting $D^{\prime}.$ Then we conclude that $f$ satisfies tangential CR equation
on $\Pi \cap \partial D.$ Due to the arbitrariness of the section $\Pi,$ the function $f$ satisfies tangential CR equation everywhere on $\partial D$ and hence is the boundary value of a holomorphic function in $D.$
\qed

The following  consequence of Corollary \ref{C:V^p} is another version of the conjecture from \cite{GS}.
\begin{Theorem} \label{T:G_S} Let $D \subset \mathbb C^n, n \geq 2,$ be a bounded domain with real analytic strictly convex boundary. Let $S \subset\subset D$ be a
real-analytic closed hypersurface. Suppose that there is an open set $U \subset D$ such that no complex line tangent to $S$ intersects
$U$ (for instance, $S$ is convex). Then any  function  $f \in C^{\omega}(\partial D),$ admitting analytic extension into each complex line tangent to $S,$
is the boundary value of a function holomorphic in $D.$
\end{Theorem}
\pf
Since $\partial D$ is strictly convex, the intersections $D \cap L$ with any complex line $L \in T^{\mathbb C}S,$ tangent to $S,$ is an analytic disc.
The family of those discs is parametrized by the complex tangent bundle of $S$ and is real-analytic and regular. The boundaries
$L \cap \partial D$ cover the whole $\partial D$ because $S$ is a closed hypersurface.

All we need to check now is that there are 3-dimensional degenerate and homologically nontrivial (having nontrivial orbits)  subfamilies passing through each point
$b \in \partial D,$ in any prescribed admissible 3-dimensional direction $\Pi.$

To this end, take a two-dimensional
complex plane $P^2$ containing $b,$ transversally intersecting $\partial D$  and such that $P^2 \cap U \neq \emptyset$ and
the tangent plane of the intersection $P^2 \cap \partial D$ is  $$T(P^2 \cap \partial D)=\Pi.$$

The intersection $P \cap S$ is 3-dimensional. For each $t \in S$ there is a unique complex line, $L_t,$ contained in the tangent plane $T_t S.$  This line
$L_t$ belongs to the complex plane $P^2.$
The mapping $$P^2\cap S \ni t \to D_t=L_t \cap D$$ defines 3-dimensional subfamily of attached analytic discs, parametrized by the hypersurface $S.$
The parametrization $\Phi$ can be chosen so that $\Phi(0,t)=t.$
\footnote{Existence of a parametrization of the form $\Phi(\zeta,t)$  for the family of complex tangent lines
requires triviality of the complex tangent bundle $T^{\mathbb C}S.$ It is so ,for instance, if $S$ is a sphere in $\mathbb C^2$.
If this is not the case then we cut the manifold $S$ into pieces with the trivial complex tangent bundles and prove the claim separately for each
corresponding bordered portion of $\partial D.$}
This subfamily is degenerate
because the boundaries $\partial D_t$ sweep up the manifold of the same dimension $d=3$ as the dimension $k=3$ of the subfamily, so we are in the situation $k+1>d$ (the case b) in Proposition \ref{P:condition_1}).

Moreover, the subfamily has the homologically nontrivial orbit because the 3-cycle $\Phi(\{0\} \times (P^2 \cap S))=P^2 \cap S$ is not homological to zero in the union
$\cup_{t \in P^2 \cap S} \overline D_t.$
Indeed,  any 4-chain $V$ bounded by $S$ must contain the "hole", $U$. But the complex lines $L_t$ do not intersect $U$ and hence $U$ belongs to the complement
of the discs $\overline D_t.$

Thus, all the conditions of Corollary \ref{C:V^p} are fullfiled and therefore $f$ extends holomorphically inside the domain $D$.
\qed
\begin{Remark}
The condition that the complex lines are tangent to the hypersurface $S$ is not essential. It is just one of the ways to define a family of attached analytic discs and gives a concrete parametrization of it. If fact, one can consider other  families, parametrized in a different way. The condition for the complex lines not to meet an open set $U$ is also one of the possible ways to impose the condition of the homological nontriviality.
\end{Remark}

\section{Concluding remarks.}

\begin{itemize}
\item Everywhere in the article the parametrizing manifold $M$ was assumed closed, i.e. having the empty boundary. Nevertheless,  the main results of the article remain true for the case of nonempty boundary as well. In this case, the condition of homological nontriviality formulates in terms of relative homologies.  However, dealing with relative homology groups and taking care about the contribution of the boundary $\partial M$ imply technical complication. Therefore, we deicided to restrict oursleves here only with the case $\partial M=\emptyset,$ which we think is  enough for the demonstration of the main ideas of the article.

\item  Corollary \ref{C:n-dim-strip} ($n$-dimensional strip-problem), for the special case of $SU(n)-$ invariant families of analytic discs in $\mathbb C^n$, was obtained in \cite{AS}.
\item
The result similar to Theorem \ref{T:Thm7} but for special families of attached analytic discs (so called thin discs)
was obtained by Tumanov
\cite{T1}. In his result the family of discs is constructed individually for each manifold $\Omega$ and depends on $\Omega.$
\item
In a recent article \cite{A_merom} a characterization of polyanalytic functions, similar to  Theorem \ref{C:4.2},  was proved.
\item
In Theorem \ref{T:Main}, we consider the $(k+2)$-dimensional manifold $\Sigma=\Delta \times M^k$ which carries a natural CR- structure.
Essentially, what we prove in Theorem \ref{T:Main} is that if a CR-mapping $\Phi$ is not degenerate (has the maximal rank) on
$\Sigma$   and if $\Phi$ induces nontrivial   homomorphism of the  homology groups $H_k(\overline \Sigma) \cong H_k(M^k)$ to
the group $H_k(\Phi(\overline \Sigma))$
(i.e. $\Phi(\zeta,t)$ has the homologically nontrivial $t-$ orbit), then for the image of the boundary $b\Sigma=S^1 \times M$
we have $H_{k+1}(\Phi(b \Sigma)) \neq 0.$
It might be interesting to generalize that statement to
CR-manifolds of  more general form (not necessarily direct products of complex disc and a real manifold).
The expected theorem might state that, under certain conditions,
nondegenerate CR-mapping $\Phi$ of a CR-manifold can not be homologically trivial on the boundary.
Moreover, analogously to the classical case (Proposition \ref{P:argument_principle},) where the statement about collapse of the interior is a corollary of the argument principle (which, in turn, can be rephrased as a rule of computing the linking numbers between two complex manifolds), the expected result might be a particular case of a more general theorem about linking numbers between two CR-manifolds.
\end{itemize}
\medskip
{\it Ackhowledgements.} The author thanks  T.-C. Dinh, P. Dolbeault, V. Gichev, J. Globevnik,
G. Henkin, S. Shnider, E. L. Stout and A. Tumanov  for useful discussions and remarks. The author is especially grateful to Peter Kuchment for his valuable editorial remarks and advices. The author thanks the referees for many useful remarks and suggestions which helped to improve the revised version of the article.
%

This work was partially supported by ISF (Israel Science Foundation), Grant 688/08.
Some of this research was done as a part of European Networking Program HCAA.

\noindent
{\it address:} Department of Mathematics, Bar-Ilan University, Ramat-Gan, 52900, Israel.

\noindent
{\it e-mail:} agranovs@macs.biu.ac.il

\end{document}